% Dmitri Pavlov's Plain TeX macros
% Copyright 2017, 2018 Dmitri Pavlov

% Remove \outer
\edef\inewcount{\noexpand\csname newcount\endcsname}
\edef\inewdimen{\noexpand\csname newdimen\endcsname}
\edef\inewskip{\noexpand\csname newskip\endcsname}
\edef\inewmuskip{\noexpand\csname newmuskip\endcsname}
\edef\inewbox{\noexpand\csname newbox\endcsname}
\edef\inewhelp{\noexpand\csname newhelp\endcsname}
\edef\inewtoks{\noexpand\csname newtoks\endcsname}
\edef\inewread{\noexpand\csname newread\endcsname}
\edef\inewwrite{\noexpand\csname newwrite\endcsname}
\edef\inewfam{\noexpand\csname newfam\endcsname}
\edef\inewlanguage{\noexpand\csname newlanguage\endcsname}
\edef\inewinsert{\noexpand\csname newinsert\endcsname}
\edef\inewif{\noexpand\csname newif\endcsname}

% UTF-8 encoding for Plain TeX.
% Copyright 2008, 2015, 2018 Dmitri Pavlov.
% You may redistribute this file under the terms of GNU General Public License version 3.
% Report bugs and suggestions to me by email (host: math.berkeley.edu, user: pavlov).
% So far this file includes support for most of TeX symbols and Latin 1 letters.

% Setup catcodes for UTF-8 characters: legitimate initial octets are active, intermediate octets are letters, other octets are illegal.
\countdef\ch=253
\ch="80 \loop\ifnum\ch<"100 \lccode\ch=0 \uccode\ch=0 \advance\ch1 \repeat % for LaTeX
\ch="80 \loop\ifnum\ch<"C0 \catcode\ch=11 \advance\ch1 \repeat % make all intermediate UTF-8 octets (top 2 bits are 10) letters
\ch="C0 \loop\ifnum\ch<"100 \catcode\ch=15 \advance\ch1 \repeat % make all octets with the top 2 bits set illegal
\ch="C2 \loop\ifnum\ch<"F5 \catcode\ch=\active \advance\ch1 \repeat % make all legitimate initial UTF-8 octets active
\ch="C2 \loop\ifnum\ch<"E2 \uccode\ch=1 \advance\ch1 \repeat % make uccode=1 for all initial UTF-8 octets of letters, to allow macros to recognize Unicode letters
\let\ch\underfined

% \cdef\macro+ is like \def\macro{+}, but the catcode of + will be 11
\catcode0=12 \def\cdef#1#2{\begingroup\lccode0=`#2 \lowercase{\endgroup \def#1{^^@}}} \catcode0=9

% Preliminary macros for UTF-8: every legitimate UTF-8 sequence produces a control sequence whose name consists precisely of the octets in the sequence.
\catcode`@=\active
\def@#1{\cdef\chr#1 \edef#1##1{\noexpand\csname\chr##1\noexpand\endcsname}} % 2-byte
@^^c2@^^c3@^^c4@^^c5@^^c6@^^c7@^^c8@^^c9@^^ca@^^cb@^^cc@^^cd@^^ce@^^cf@^^d0@^^d1@^^d2@^^d3@^^d4@^^d5@^^d6@^^d7@^^d8@^^d9@^^da@^^db@^^dc@^^dd@^^de@^^df
\def@#1{\cdef\chr#1 \edef#1##1##2{\noexpand\csname\chr##1##2\noexpand\endcsname}} % 3-byte
@^^e0@^^e1@^^e2@^^e3@^^e4@^^e5@^^e6@^^e7@^^e8@^^e9@^^ea@^^eb@^^ec@^^ed@^^ee@^^ef
\def@#1{\cdef\chr#1 \edef#1##1##2##3{\noexpand\csname\chr##1##2##3\noexpand\endcsname}} % 4-byte
@^^f0@^^f1@^^f2@^^f3@^^f4
\let\chr\undefined
\let\cdef\undefined

% The following three lines establish a concise syntax for defining characters: each line contains a UTF-8 character followed by its definition.
\def\grabfuturelet{\futurelet\next\grabexamine}
\def\grabexamine{\ifx\next\csname\expandafter\grab\fi}
\obeylines \def\grab\csname#1\endcsname#2^^M{\expandafter\def\csname#1\endcsname{#2}\expandafter\grabfuturelet} \expandafter\grabfuturelet%
% Coverage: all TeX symbols except those not present in Unicode and the part of Latin 1 and Latin Extended-A that is present in the Computer Modern fonts.
 ~
¢{\hbox{\rm\rlap/c}}
£{\it\$}
%¤
%¥
%¦
%¨
%ª
«\leftguillemet
­\-
%®
%¯
%°
%²
%³
%´
%µ
%·\cdotp
%¸\c
%¹
%º
»\rightguillemet
%¼
%½
%¾
À{\`A}
Á{\'A}
Â{\^A}
Ã{\~A}
Ä{\"A}
Ç{\c C}
È{\`E}
É{\'E}
Ê{\^E}
Ë{\"E}
Ì{\`I}
Í{\'I}
Î{\^I}
Ï{\"I}
% must define \ETH, \THORN, \eth, \thorn separately
Ð\ETH
Ñ{\~N}
Ò{\`O}
Ó{\'O}
Ô{\^O}
Õ{\~O}
Ö{\"O}
Ù{\`U}
Ú{\'U}
Û{\^U}
Ü{\"U}
Ý{\'Y}
Þ\THORN
à{\`a}
á{\'a}
â{\^a}
ã{\~a}
ä{\"a}
ç{\c c}
è{\`e}
é{\'e}
ê{\^e}
ë{\"e}
ì{\`\i}
í{\'\i}
î{\^\i}
ï{\"\i}
ð\eth
ñ{\~n}
ò{\`o}
ó{\'o}
ô{\^o}
õ{\~o}
ö{\"o}
ù{\`u}
ú{\'u}
û{\^u}
ü{\"u}
ý{\'y}
þ\thorn
ÿ{\"y}
% Latin Extended-A
Ā{\=A}
ā{\=a}
Ă{\u A}
ă{\u a}
%Ą
%ą
Ć{\'C}
ć{\'c}
Ĉ{\^C}
ĉ{\^c}
Ċ{\.C}
ċ{\.c}
Č{\v C}
č{\v c}
Ď{\v D}
ď{\v d}
%Đ
%đ
Ē{\=E}
ē{\=e}
Ĕ{\u E}
ĕ{\u e}
Ė{\.E}
ė{\.e}
%Ę
%ę
Ě{\v E}
ě{\v e}
Ĝ{\^G}
ĝ{\^g}
Ğ{\u G}
ğ{\u g}
Ġ{\.G}
ġ{\.g}
Ģ{\c G}
ģ{\c g}
Ĥ{\^H}
ĥ{\^h}
%Ħ
%ħ
Ĩ{\~I}
ĩ{\~\i}
Ī{\=I}
ī{\=\i}
Ĭ{\u I}
ĭ{\u\i}
%Į
%į
İ{\.I}
%ı\i
ĲIJ
ĳij
Ĵ{\^J}
ĵ{\^\j}
Ķ{\c K}
ķ{\c k}
%ĸ
Ĺ{\'L}
ĺ{\'l}
Ļ{\c L}
ļ{\c l}
Ľ{\v L}
ľ{\v l}
%Ŀ
%ŀ
%Ł\L
%ł\l
Ń{\'N}
ń{\'n}
Ņ{\c N}
ņ{\c n}
Ň{\v N}
ň{\v n}
%ŉ
%Ŋ
%ŋ
Ō{\=O}
ō{\=o}
Ŏ{\u O}
ŏ{\u o}
Ő{\"O}
ő{\"o}
%Œ\OE
%œ\oe
Ŕ{\'R}
ŕ{\'r}
Ŗ{\c R}
ŗ{\c r}
Ř{\v R}
ř{\v r}
Ś{\'S}
ś{\'s}
Ŝ{\^S}
ŝ{\^s}
Ş{\c S}
ş{\c s}
Š{\v S}
š{\v s}
Ţ{\c T}
ţ{\c t}
Ť{\v T}
ť{\v t}
%Ŧ
%ŧ
Ũ{\~U}
ũ{\~u}
Ū{\=U}
ū{\=u}
Ŭ{\u U}
ŭ{\u u}
%Ů
%ů
Ű{\H U}
ű{\H u}
%Ų
%ų
Ŵ{\^W}
ŵ{\^w}
Ŷ{\^Y}
ŷ{\^y}
Ÿ{\"Y}
Ź{\'Z}
ź{\'z}
Ż{\.Z}
ż{\.z}
Ž{\v Z}
ž{\v z}
%ſ
’'
‘`
”{''}
“{``}
‐-
–{--}
—{---}
% these ligatures should not be used in TeX text
%ﬀ{ff}
%ﬁ{fi}
%ﬂ{fl}
%ﬃ{ffi}
%ﬄ{ffl}
¡{!`}
¿{?`}
−-
′'
% these ligatures should not be used in TeX text
%″{''}
%‴{'''}
%⁗{''''}
ß\ss
æ\ae
Æ\AE
œ\oe
Œ\OE
ø\o
Ø\O
å\aa
Å\AA
ł\l
Ł\L
% defined below
%ı\i
%ȷ\j
†\dag
‡\ddag
§\S
¶\P
% combining characters are left undefined because they come after, not before the accented character
%◌̣\d
%◌̱\b
%◌̧\c
©\copyright
…\dots
%◌̀\`
%◌́\'
%◌̌\v
%◌̆\u
%◌̄\=
%◌̂\^
%◌̇\.
%◌̋\H
%◌̃\~
%◌̈\"
%◌͡\t
%\brace[lr][du]
α\alpha
β\beta
γ\gamma
δ\delta
ϵ\epsilon
ζ\zeta
η\eta
θ\theta
ι\iota
κ\kappa
λ\lambda
μ\mu
ν\nu
ξ\xi
οo
π\pi
ρ\rho
σ\sigma
τ\tau
υ\upsilon
ϕ\phi
χ\chi
ψ\psi
ω\omega
ε\varepsilon
ϑ\vartheta
ϖ\varpi
ϱ\varrho
ς\varsigma
φ\varphi
Γ\Gamma
Δ\Delta
Θ\Theta
Λ\Lambda
Ξ\Xi
Π\Pi
Σ\Sigma
Υ\Upsilon
Φ\Phi
Ψ\Psi
Ω\Omega
ℵ\aleph
ℏ\hbar
ı\relax\ifmmode\imath\else\i\fi
ȷ\relax\ifmmode\jmath\else\j\fi
ℓ\ell
℘\wp
ℜ\Re
ℑ\Im
∂\partial
∞\infty
%'\prime
∅\emptyset
∇\nabla
√\surd
⊤\top
⊥\bot
∠\angle
△\triangle
∀\forall
∃\exists
¬\neg
♭\flat
♮\natural
♯\sharp
♣\clubsuit
♢\diamondsuit
♡\heartsuit
♠\spadesuit
∐\coprod
⋁\bigvee
⋀\bigwedge
⨄\biguplus
⋂\bigcap
⋃\bigcup
∫\int
%!\intop
∏\prod
∑\sum
⨂\bigotimes
⨁\bigoplus
⨀\bigodot
∮\oint
%!\ointop
⨆\bigsqcup
%!\smallint
◁\triangleleft
▷\triangleright
△\bigtriangleup
▽\bigtriangledown
∧\wedge
∨\vee
∩\cap
∪\cup
%‡\ddagger : \ddag
%†\dagger : \dag
⊓\sqcap
⊔\sqcup
⊎\uplus
⨿\amalg
⋄\diamond
∙\bullet
≀\wr
÷\div
⊙\odot
⊘\oslash
⊗\otimes
⊖\ominus
⊕\oplus
∓\mp
±\pm
∘\circ
%white circle:
○\Orb
%large circle:
◯\bigcirc
∖\setminus
⋅\cdot
∗\ast
% Latin 1
×\times
% Unicode
⨯\times
⋆\star
∝\propto
⊑\sqsubseteq
⊒\sqsupseteq
∥\parallel
‖\|
∣\divides
%|\mid
⊣\dashv
⊢\vdash
↗\nearrow
↘\searrow
↖\nwarrow
↙\swarrow
⇔\Leftrightarrow
⇐\Leftarrow
⇒\Rightarrow
≠\neq
≤\leq
≥\geq
≻\succ
≺\prec
≈\approx
≽\succeq
≼\preceq
⊃\supset
⊂\subset
⊇\supseteq
⊆\subseteq
∈\in
∋\ni
≫\gg
≪\ll
%◌̸\not
↔\leftrightarrow
←\leftarrow
→\rightarrow
↦\mapsto
%↦\mapstochar : \mapsto
∼\sim
≃\simeq
%⊥\perp : \bot
≡\equiv
≍\asymp
⌣\smile
⌢\frown
↼\leftharpoonup
↽\leftharpoondown
⇀\rightharpoonup
⇁\rightharpoondown
↪\hookrightarrow
↩\hookleftarrow
%\lhook
%\rhook
⋈\bowtie
⊨\models
⟹\Longrightarrow
⟶\longrightarrow
⟵\longleftarrow
⟸\Longleftarrow
⟼\longmapsto
⟷\longleftrightarrow
⟺\Longleftrightarrow
%⇔\iff : \Longleftrightarrow
%.\ldotp : .
%⋅\cdotp : \cdot
%:\colon : :
%…\ldots : \dots
⋯\cdots
⋮\vdots
⋱\ddots
%\acute, \grave, \ddot, \tilde, \bar, \breve, \check, \hat, \vec, \dot, \widetilde, \widehat
%\overrightarrow, \overleftarrow, \overbrace, \underbrace
%\lmoustache, \rmoustache, \lgroup, \rgroup, \arrowvert, \Arrowvert, \bracevert
∥\Vert
%|\vert
↑\uparrow
↓\downarrow
↕\updownarrow
⇑\Uparrow
⇓\Downarrow
⇕\Updownarrow
%\\backslash
⟩\rangle
⟨\langle
%{\lbrace
%}\rbrace
⌉\rceil
⌈\lceil
⌋\rfloor
⌊\lfloor
≅\cong
∉\notin
⇌\rightleftharpoons
≐\doteq
% Combined symbols from Unicode math blocks
∄\not\exists
∌\not\ni
∔\dot+
∕/
∣|
∤\not|
∦\not\|
∬\int\!\!\!\int
∭\int\!\!\!\int\!\!\!\int
∮\oint
∸\dot-
≁\not\sim
≄\not\simeq
≆\not\cong
≇\not\cong
≉\not\approx
≐\dot=
≔:=
≕=:
≢\not\equiv
≭\not\asump
≮\not<
≯\not<
≰\not\le
≱\not\ge
⊀\not\prec
⊁\not\succ
⊄\not\subset
⊅\not\supset
⊈\not\subseteq
⊉\not\supseteq
⊦\vdash
⊧\models
⊬\not\vdash
⊭\not\models
⊲\triangleleft
⊳\triangleright
⋠\not\preceq
⋡\not\succeq
⋤\not\sqsubseteq
⋥\not\sqsupseteq
⋪\not\triangleleft
⋫\not\triangleright
◻\square

\catcode`\^^M=5 %
\let\grabfuturelet\undefined \let\grabexamine\undefined \let\grab\undefined

% Generic macros for UTF-8: every legitimate UTF-8 sequence produces a control sequence whose name consists precisely of the octets in the sequence.
% An undefined control sequence produces an error message.
\let\xcsname=\csname
\let\xendcsname=\endcsname
\def@#1{\def#1##1{\expandafter\ifx\csname\string#1##1\endcsname\relax\errmessage{Undefined UTF-8 sequence \string#1##1}\else\xcsname\string#1##1\xendcsname\fi}}
@^^c2@^^c3@^^c4@^^c5@^^c6@^^c7@^^c8@^^c9@^^ca@^^cb@^^cc@^^cd@^^ce@^^cf@^^d0@^^d1@^^d2@^^d3@^^d4@^^d5@^^d6@^^d7@^^d8@^^d9@^^da@^^db@^^dc@^^dd@^^de@^^df
\def@#1{\def#1##1##2{\expandafter\ifx\csname\string#1##1##2\endcsname\relax\errmessage{Undefined UTF-8 sequence \string#1##1##2}\else\xcsname\string#1##1##2\xendcsname\fi}}
@^^e0@^^e1@^^e2@^^e3@^^e4@^^e5@^^e6@^^e7@^^e8@^^e9@^^ea@^^eb@^^ec@^^ed@^^ee@^^ef
\def@#1{\def#1##1##2##3{\expandafter\ifx\csname\string#1##1##2##3\endcsname\relax\errmessage{Undefined UTF-8 sequence \string#1##1##2##3}\else\xcsname\string#1##1##2##3\xendcsname\fi}}
@^^f0@^^f1@^^f2@^^f3@^^f4
\let@\undefined
\catcode`@=12

\newif\ifscroll % scroll vs codex
\newif\ifsuppressunusedbib % suppress unused bibliography items?

% Diagonostics
\tracinglostchars=2

\def\printerr#1{\immediate\write17{#1}}
\def\warningline#1#2{\printerr{! #2}\printerr{l.#1}\printerr{}}
\def\ewarningline#1#2#3{\printerr{! #2}\printerr{l.#1 #3}\printerr{}}
\def\warning{\warningline{\the\inputlineno}}
%\let\warning\errmessage % turn warnings into errors

% Syntactic macros
\long\def\gobble#1{}
\ifx\gobbleinit\undefined{\long\gdef\gobbleinit#1\par{}}\fi
\def\expand#1{\edef\expandmacro{#1}\expandmacro\let\expandmacro\undefined}
\def\setetok#1#2{\expand{\noexpand#1{#2}}}
\def\expandtoks#1{\expandafter\edef\expandafter\expandmacro\expandafter{\the#1}#1\expandafter{\expandmacro}}
\def\appendexpand#1#2{\setetok#1{\the#1#2}}
\long\def\append#1#2{#1\expandafter{\the#1#2}}
\long\def\appendtoksexpand#1#2{#1\expandafter\expandafter\expandafter{\expandafter\the\expandafter#1\the#2}}
\long\def\appendonceexpand#1#2{#1\expandafter\expandafter\expandafter{\expandafter\the\expandafter#1#2}}

% Macros that use \special, known to work with xdvi, dvips, dvipdfm

%\let\printlink\print % print links and anchors
\def\link#1#2{\lhighlight{#2}}
\def\llink#1{\printlink{llink #1}\link{\ohash#1}}
\catcode`\#=11 \def\ohash{#}\catcode`\#=6
\catcode`\&=11 \def\ampersand{&}\catcode`\&=4
\def\anchor#1#2{\printlink{anchor #1 #2}#2}

\def\setpapersize#1#2{} % \number#1 sp triggers bugs in dvips and does not work with dvipdfm(x)
\def\dumpbox#1#2#3{\shipout\vbox{\setpapersize{#1}{#2}\unvbox#3}}
\def\mps#1{\epsfbox{#1}}
\def\metadata#1#2{}
\def\src{} % for proper match between bibliographic references and source code

% A simplified version of the epsfbox macro from epsf.tex, exclusively for METAPOST output, with HiResBoundingBox
% Normally should be invoked through \mps, which works also with PDF
\newread\epsffilein
\newif\ifepsfbbfound\inewif\ifepsffilecont
\newdimen\epsfxsize\inewdimen\epsfysize
\newdimen\pspoints\pspoints1bp
\let\runmp\errmessage % will be set to \warning if the METAPOST file cannot be compiled
\def\epsfbox#1{\openin\epsffilein=#1 \ifeof\epsffilein\runmp{Could not open file #1}\else
	{\def\do##1{\catcode`##1=12}\dospecials\catcode`\ =10\epsffileconttrue
		\epsfbbfoundfalse
		\loop\read\epsffilein to\epsffileline \ifeof\epsffilein\epsffilecontfalse\else\expandafter\epsfaux\epsffileline :. \\\fi\ifepsffilecont\repeat
		\ifepsfbbfound\else\errmessage{No HiResBoundingBox comment found in file #1}\fi}%
	\closein\epsffilein
	\epsfysize\epsfury\pspoints \advance\epsfysize-\epsflly\pspoints
	\epsfxsize\epsfurx\pspoints \advance\epsfxsize-\epsfllx\pspoints
	% We create a box with height and depth corresponding to the two vertical dimensions in the bounding box and width given by the total width
	\setbox0\hbox{\vbox to\epsfury\pspoints{\vfil\hbox to\epsfxsize{\dimen0=\epsfllx\pspoints \kern-\dimen0 \includegraphics{#1}\hfil}}}%
	\dp0=\epsflly\pspoints \dp0=-\dp0
	\box0 \fi}
%	\hbox{\vbox to\epsfysize{\vfil\hbox to\epsfxsize{\special{psfile=#1}\hfil}}}}
% llx=\epsfllx\space lly=\epsflly\space urx=\epsfurx\space ury=\epsfury
%	    rwi=\number\epsftmp\space rhi=\number\epsfrsize
\catcode`\%=12 \def\epsfbblit{%%HiResBoundingBox} \catcode`\%=14
\let\dummybrace=} % to compensate for { on the previous line when it is read by \read
\def\epsfaux#1:#2\\{\def\testit{#1}\ifx\testit\epsfbblit \epsfgrab #2 . . . \\\epsffilecontfalse\epsfbbfoundtrue\fi}
\def\empty{}
\def\epsfgrab #1 #2 #3 #4 #5\\{\gdef\epsfllx{#1}\ifx\epsfllx\empty\epsfgrab #2 #3 #4 #5 .\\\else\gdef\epsflly{#2}\gdef\epsfurx{#3}\gdef\epsfury{#4}\fi} % gdef because epxfbox wraps us in a group

% pdfTeX completely ignores established conventions for \special
\newif\ifpdf \pdffalse \ifx\pdfoutput\undefined\else\ifx\pdfoutput\relax\else\ifnum\pdfoutput<1 \else\pdftrue\fi\fi\fi
\ifpdf
\pdfcompresslevel=0 % HTTP offers better compression methods anyway
\pdfobjcompresslevel=0
\def\pdflink#1#2{\leavevmode \lhighlight{\pdfstartlink user { /Subtype /Link /Border [0 0 0] /A << /S #1 >> }#2\pdfendlink}}
\def\llink#1{\pdflink{/GoTo /D (#1)}}
\def\link#1{\pdflink{/URI /URI (#1)}}
\def\anchor#1#2{\pdfdest name {#1} xyz #2}

\def\setpapersize#1#2{\pdfpagewidth#1 \pdfpageheight#2 }
\def\dumpbox#1#2#3{\setpapersize{#1}{#2}\shipout\box#3}
\def\metadata#1#2{\pdfinfo{/Title (#1) /Author (#2)}}
\input supp-pdf % the only external dependence of DPMAC: a converter from METAPOST's output to PDF
\def\mps#1{\convertMPtoPDF{#1}{1}{1}}
%\chardef\makeMPintoPDFobject=1
\fi

% Embedded METAPOST
\newtoks\buffertoks
\def\addcode{\immediate\write\mpout}
\def\addunexpandedcode#1{{\toks0={#1}\addcode{\the\toks0}}}
\def\addcodebuffer#1{\edef\tmp{#1}\buffertoks\expandafter\expandafter\expandafter{\expandafter\the\expandafter\buffertoks\tmp}}
\def\addunexpandedcodebuffer#1{\buffertoks\expandafter{\the\buffertoks#1}}

\def\grabcode{\catcode`\#=12 \endlinechar=10 
	\afterassignment\dumpcode\outtoks} % make sure the newline characters get written
\def\dumpcode{\addcode{\the\outtoks}\endinlinemp\gobble} % \gobble the newline with character code 10
\def\begininlinemp{\inimp\begingroup\catcode`\^=7 \iftypesetting\mps{\filestem.\the\figno}\let\addcode\gobble\fi \addcode{beginfig(\the\figno);}}
\def\endinlinemp{\addcode{endfig;}\addcode{}\endgroup\global\advance\figno1\relax}
\ifx\endprolog\undefined\let\endprolog\relax\fi

% Repair LaTeX damage
\def\plainfmtname{plain}\ifx\fmtname\plainfmtname\else
\edef\plainoutput{\the\output}
\global\chardef\itfam=4
\def\_{\leavevmode \kern.06em \vbox{\hrule width.3em}} % for LaTeX

\outputpenalty=0
\tracingstats=0
\newlinechar=-1
\maxdeadcycles=25
\showboxbreadth=5
\showboxdepth=3
\errorcontextlines=5
\overfullrule=5pt
\maxdepth=4pt
\parindent=20pt
\abovedisplayskip=12pt plus 3pt minus 9pt
\belowdisplayskip=12pt plus 3pt minus 9pt
\belowdisplayshortskip=7pt plus 3pt minus 4pt

 % math italic

\font\tensy=cmsy10

\font\sevenrm=cmr7

\catcode"18=12
\catcode`@=11
{
\global\let\end\@@end
\global\let\input\@@input
}
%\let\everymath\frozen@everymath
%\let\everydisplay\frozen@everydisplay
% copied from plain.tex:
\def\eqalign#1{\null\,\vcenter{\openup\jot\m@th
  \ialign{\strut\hfil$\displaystyle{##}$&$\displaystyle{{}##}$\hfil
      \crcr#1\crcr}}\,}
\catcode`@=12
%\def\line{\hbox to\hsize}
%\everymath{}
%\everydisplay{}
\fi

% AMS fonts: family 8, 9, 10
\font\tenmsa=msam10 \font\sevenmsa=msam7 \font\fivemsa=msam5 \newfam\msafam \textfont\msafam=\tenmsa \scriptfont\msafam=\sevenmsa \scriptscriptfont\msafam=\fivemsa % 8
   %\newfam\msbfam \textfont\msbfam=\tenmsb \scriptfont\msbfam=\sevenmsb \scriptscriptfont\msbfam=\fivemsb
%\def\Bbb{\fam\msbfam}
\font\teneufm=eufm10 \font\seveneufm=eufm7 \font\fiveeufm=eufm5 \newfam\eufmfam \textfont\eufmfam=\teneufm \scriptfont\eufmfam=\seveneufm \scriptscriptfont\eufmfam=\fiveeufm % 9

\font\teneufb=eufb10 \font\seveneufb=eufb7 \font\fiveeufb=eufb5 \newfam\eufbfam \textfont\eufbfam=\teneufb \scriptfont\eufbfam=\seveneufb \scriptscriptfont\eufbfam=\fiveeufb % 10

\font\teneurm=eurm10 \font\seveneurm=eurm7 \font\fiveeurm=eurm5 \newfam\eurmfam \textfont\eurmfam=\teneurm \scriptfont\eurmfam=\seveneurm \scriptscriptfont\eurmfam=\fiveeurm % 11

\font\teneurb=eurb10 \font\seveneurb=eurb7 \font\fiveeurb=eurb5 \newfam\eurbfam \textfont\eurbfam=\teneurb \scriptfont\eurbfam=\seveneurb \scriptscriptfont\eurbfam=\fiveeurb % 12

\font\teneusm=eusm10 \font\seveneusm=eusm7 \font\fiveeusm=eusm5 \newfam\eusmfam \textfont\eusmfam=\teneusm \scriptfont\eusmfam=\seveneusm \scriptscriptfont\eusmfam=\fiveeusm % 13

\font\teneusb=eusb10 \font\seveneusb=eusb7 \font\fiveeusb=eusb5 \newfam\eusbfam \textfont\eusbfam=\teneusb \scriptfont\eusbfam=\seveneusb \scriptscriptfont\eusbfam=\fiveeusb % 14
\def\eucalbf{\fam\eusbfam}
%\font\teneuex=euex10 \font\seveneuex=euex7 \newfam\euexfam \textfont\euexfam=\teneuex \scriptfont\euexfam=\seveneuex

\font\tenss=cmss10 \font\sevenss=cmss7 \font\fivess=cmss5 \newfam\ssfam \textfont\ssfam\tenss \scriptfont\ssfam\sevenss \scriptscriptfont\ssfam\fivess % 15
\def\sf{\fam\ssfam}

\font\sevenit=cmti7 \scriptfont\itfam=\sevenit

% Font style
\let\articletitle\seventeenss
\let\chaptertitle\twelvebf
\let\sectiontitle\tenbf
\let\subsectiontitle\tenbfit
\let\subsubsectiontitle\tenit
\let\contchaptertitle\tenbf % chapter titles in the table of contents
\let\contsectiontitle\tenrm % ditto for sections
\let\contsubsectiontitle\sevenrm % ditto for subsections
\let\contsubsubsectiontitle\fiverm % ditto for subsubsections
\let\parnumfont\tenrm % paragraph numbers
\let\parbackreffont\fiverm % back references at the end of paragraphs
\let\proclaimfont\tenbf
\let\prooffont\tenit
\let\mainfont\tenrm

% Formatting of draft comments

\def\lhighlight{} % in principle, links can be highlighted, but this is usually done by the renderer
 % for screen
 % for printing drafts
\def\suppresscomments#1#2#3{} % for sharing drafts
\def\prohibitcomments#1#2#3{\errmessage{Draft comments are not allowed in the final version}} % for final versions

\ifx\format\undefined\else\format\fi % default paper size is letter
\ifx\comment\undefined\def\comment{\prohibitcomments}\fi

\newdimen\hmargin
\hmargin=1in
\newdimen\vmargin
\vmargin=1in
\newdimen\plaintextwidth
\plaintextwidth=6.5in
\newdimen\plaintextheight
\plaintextheight=8.9in

% Output routine for scrolls with variable height
\newdimen\totalht
\newdimen\totalwd
\def\shipbox#1{%
	\totalht\ht#1
	\advance\totalht2\vmargin
	\totalwd\wd#1
	\advance\totalwd2\hmargin
	\hoffset-1in
	\advance\hoffset\hmargin
	\voffset-1in
	\advance\voffset\vmargin
	\dumpbox\totalwd\totalht{#1}}

% Contents
\newtoks\cont
\def\contents{\begingroup\suppressbackreftrue\the\cont\endgroup}
\let\printcont\gobble
\def\contlinechapt#1#2{\printcont{#2}\smallskip\everypar{}\noindent{\contchaptertitle\llink{chapter.#1}{#2}}\par} % Link: contents to chapter
\def\contlinesect#1#2{\printcont{#1.  #2}\everypar{}\indent{\contsectiontitle\llap{#1.\enskip}\llink{section.#1}{#2}}\par} % Link: contents to section
\def\contlinesubsect#1#2{\printcont{#1.  #2}\everypar{}\indent{\contsubsectiontitle{#1.\enskip}\llink{paragraph.#1}{#2}}\par} % Link: contents to subsection
\def\contlinesubsubsect#1#2{\printcont{#1.  #2}\everypar{}\indent\indent\indent{\contsubsubsectiontitle\llap{#1.\enskip}\llink{paragraph.#1}{#2}}\par} % Link: contents to subsubsection
\def\addcont#1#2#3{\append\cont{#1}\appendexpand\cont{{#2}}\append\cont{{#3}}}

% Chapters
\newif\ifpresec \presecfalse % true if we are in a chapter, but before its first section
\def\chapter#1\par{%
	\def\chapname{#1}%
	\parn0
	\subparn0
	\presectrue
	\numbfalse
	\addcont\contlinechapt\chapname{#1}%
	\curverb{Chapter~}%
	\assignlabel{chapter}{\chapname}%
	\tchapter#1\par}
\def\tchapter#1\par{% typeset chapter title
	\everypar{\let\beforesect\beforesection}% do not number paragraphs before the first section
	\chapbreak\bigbreak
	\centerline{\plabel\chaptertitle\anchor{chapter.#1}{#1}}% Anchor: chapter
	\nobreak\medskip
	\let\beforesect\relax % do not break pages right after a chapter title
}

\let\chapbreak\relax % \let\chapbreak\tchapbreak to start chapters on a new page

% Sections
\newcount\secn \secn0
\def\section#1\par{%
	\ifx\sectionid\undefined\advance\secn1 \edef\sectionid{\the\secn}\fi
	\everypar{\numpar}% number paragraphs inside a section
	\parn0
	\subparn0
	\presecfalse
	\numbfalse
	\addcont\contlinesect\sectionid{#1}%
	\curverb{\S}%
	\assignlabel{section}{\sectionid}%
	\tsection#1\par
}
\def\tsection#1\par{
	\beforesect\let\beforesect\beforesection
	\typesetsection{#1}%
	\aftersection	
	\let\sectionid\undefined
}
\def\beforesection{\vskip0pt plus.3\vsize \penalty-250 \vskip0pt plus-.3\vsize \bigskip \vskip\parskip}
\def\aftersection{\nobreak\smallskip}
\def\typesetsection#1{\leftline{\sectiontitle\ifx\sectionid\undefined\indent\else\hbox to \parindent{\hss\plabel\anchor{section.\sectionid}{\sectionid}\enspace\hfill}\fi#1}} % Anchor: section
\let\beforesect\beforesection

% Subsections
\def\subsection#1\par{\bigbreak\numbtrue\curverb{\S}\subsectiontitle\noindent#1\/\mainfont\par\nobreak\medskip\addcont\contlinesubsect{\the\secn.\the\parn}{#1}}
\def\subsubsection#1\par{\bigbreak\numbtrue\curverb{\S}\subsubsectiontitle\noindent#1\/\mainfont\par\nobreak\medskip\addcont\contlinesubsubsect{\the\secn.\the\parn}{#1}}

% Theorems and proofs
\def\sskip#1{\ifdim\lastskip<\medskipamount \removelastskip\penalty#1\medskip\fi}
\def\slug{\hbox{\kern1.5pt\vrule width2.5pt height6pt depth1.5pt\kern1.5pt}}
\newif\ifqed \newif\ifneedqed
\def\qed{\unskip\nobreak\ \slug\ifhmode\spacefactor3000 \fi\global\qedtrue}
\def\proclaim{\medbreak\atendpar{\sskip{55}}\numbtrue\gproclaim\proclaimfont} % always number theorems etc.
\def\proof{\medbreak\atendpar{\sskip{-55}}\needqedtrue\gproclaim\prooffont}
\def\xproclaim#1.{\medbreak\atendpar{\sskip{55}}{\everypar{}\noindent}{\proclaimfont#1.\enspace}\ignorespaces} % like \proclaim, no numbering
\let\abstract\xproclaim % abstracts don't get numbered

% Paragraphs
 % pseudo paragraph, cannot be labeled, does not have backrefs
 % ersatz paragraph with numbering, cannot be labeled, does not have backrefs
\def\ppar{\endgraf{\everypar{}\indent}} % paragraph inside a proof, not seen by \proof, has backrefs
\newtoks\atendpar % tokens to be inserted after \endgraf
\newtoks\atendbr % paragraph backreferences
\newif\ifnumb \numbfalse % number this paragraph?
\def\finishpar{\ifhmode\ifneedqed\ifqed\else\qed\fi\qedfalse\needqedfalse\fi\iflist\endlist\fi\the\atendbr\atendbr{}\endgraf\the\atendpar\atendpar{}\numbfalse\fi}
\def\endlist{\iflist\listfalse\endgraf{\parskip\smallskipamount\everypar{}\noindent}\fi} % terminate a list without starting a new paragraph
\newcount\parn
\newcount\subparn
\newif\ifparbref
\def\nextpar{\ifnumb\advance\parn1 \printlabel{advancing paragraph number to \the\parn}\def\brt{}%
	\ifparbref\edef\cseq{\csname backreference-list.paragraph.\the\secn.\the\parn\endcsname}%
	\expandafter\ifx\cseq\relax\else\edef\brt{\cseq}\printbackref{back references for paragraph.\the\secn.\the\parn: \cseq}\fi\fi
	\assignlabel{paragraph}{\the\secn.\the\parn}\fi}
\def\numpar{\ifnumb\nextpar\expand{\noexpand\typesetparnum{\the\secn.\the\parn}\noexpand\typesetpbr{\brt}}\fi}
\newdimen\pindent \pindent\parindent

% Paragraph and theorem numbering
\ifx\draftnum\undefined % normal numbering
\def\gproclaim#1#2.{\curverb{#2~}% #1: font, #2: prefix
	\ifnumb\nextpar\fi
	{\everypar{}\noindent}%
	\plabel#1#2%
	\ifnumb\ \anchor{paragraph.\the\secn.\the\parn}{}\the\secn.\the\parn\fi.\enspace\mainfont % Anchor: proclaim
	\ifnumb\expand{\noexpand\typesetpbr{\brt}}\fi
	\ignorespaces}
\def\typesetparnum#1{\ifnumb{\plabel\parnumfont\anchor{paragraph.#1}{}#1.\enspace}\fi} % Anchor: explicitly numbered paragraph
\def\typesetpbr#1{\ifnumb\def\brtext{#1}\ifx\brtext\empty\else\setetok\atendbr{{\parbackreffont Used in \noexpand\stripcomma\brtext.}}\fi\fi}
\else % number everything for proofreading purposes
\let\numbfalse\relax \numbtrue
\def\gproclaim#1#2.{\curverb{#2~}\medbreak\noindent#1#2.\enspace\mainfont\ignorespaces}
\inewdimen\brwidth \brwidth.6in
\parindent0pt \parskip1ex plus 1ex minus 1ex
\def\typesetparnum#1{\ifnumb\llap{\plabel\anchor{paragraph.#1}{}\parnumfont#1\enspace}\fi} % Anchor: implicitly numbered paragraph
\def\typesetpbr#1{\ifnumb\def\brtext{#1}\ifx\brtext\empty\else
	\llap{\smash{\vtop{\everypar{}\raggedright\rightskip0pt plus 0pt \leftskip0pt plus 1fill \hsize\brwidth
	\parnumfont \strut \break % the first line contains paragraph number
	\parbackreffont\stripcomma#1}}\enspace}\fi\fi}
\fi

% Lists
\def\hang{\hangindent\pindent}
\newif\iflist
\def\textindent#1{{\everypar{}\parindent\pindent\indent}\llap{#1\enspace}\listtrue\ignorespaces}

\def\li{\item{$\bullet$}}

% \par replaced by \endgraf
\def\item{\endgraf\hang\textindent}
\def\filbreak{\endgraf\vfil\penalty-200\vfilneg}

\def\eject{\endgraf\break}
\def\supereject{\endgraf\penalty-20000}
\def\smallbreak{\endgraf\ifdim\lastskip<\smallskipamount
	\removelastskip\penalty-50\smallskip\fi}
\def\medbreak{\endgraf\ifdim\lastskip<\medskipamount
	\removelastskip\penalty-100\medskip\fi}
\def\bigbreak{\endgraf\ifdim\lastskip<\bigskipamount
	\removelastskip\penalty-200\bigskip\fi}

% Email obfuscation (source and output)
%\newdimen\plht \setbox0\hbox{\char`\_} \plht\ht0 \setbox1\hbox{.} \advance\plht-\ht1 % compute the amount by which a dot accent should be lowered to become a period
%\catcode`\.\active \def.{\lower\plht\hbox{\char`\_}} \catcode`\.=12 % emails are typeset with a lowered dot accent instead of a period
%\def\email{\bgroup\catcode`.\active\xemail}
%\def\xemail#1#2{\rlap{\hphantom{#2}@#1}#2\hphantom{@#1}\egroup}

% Back references for labels and bibliography
\let\printbackref\gobble
\def\predefbackref#1{%
	\printbackref{defining back reference list backreference-list.#1}%
	\expandafter\gdef\expandafter\cseq\expandafter{\csname backreference-list.#1\endcsname}
	\expandafter\ifx\cseq\relax\expandafter\gdef\cseq{}\else\printbackref{duplicate omitted}\fi
}
\newcount\backref \backref0
\newif\ifsuppressbackref \suppressbackreffalse
\def\firstletter#1#2\endletter{#1}
\newtoks\backreflist
\newif\ifaddbr \addbrfalse
\def\recordbackref#1{%
	\edef\params{{\ifpresec\expandafter\firstletter\chapname\endletter\else\the\secn\fi.\the\parn\ifnumb\else*\fi}{\the\backref}{\the\inputlineno}}%
	\printbackref{recording back reference \string#1 for future processing with params \params}%
	\edef\tmp{\the\backreflist\noexpand\processbackref\noexpand#1\params}%
	\global\backreflist\expandafter{\tmp}%
	\printbackref{new content of backreflist: \the\backreflist}} % \global is needed because this can be invoked inside {\it ...}, say
\def\processbackref#1#2#3#4{%
	\edef\key{\expandafter\gobble\string#1}%
	\edef\cseq{\csname id.\key\endcsname}%
	\expandafter\ifx\cseq\relax \ewarningline{#4}{Undefined reference \string#1.}{\string#1}\else
		\edef\lseq{\csname backreference-list.\cseq\endcsname}%
		\printbackref{processing back reference number #3 \string#1, originating from #2, at line #4; adding to \lseq}%
		% Check for duplicates
		\edef\lastnumber{\csname lastnumber.\cseq\endcsname}%
		\edef\newnumber{#2}%
		\addbrtrue
		\printbackref{lastnumber: \lastnumber; newnumber: \newnumber;}%
		\expandafter\ifx\csname lastnumber.\cseq\endcsname\relax % we are the first back reference
		\else\ifx\lastnumber\newnumber % same as the last one
			\printbackref{Suppressing duplicate back reference #2.}%
			\addbrfalse
		\fi\fi
		\expandafter\edef\csname lastnumber.\cseq\endcsname{#2}% record the new number
		\ifaddbr
			\expandafter\expandafter\expandafter\gdef\expandafter\expandafter\csname backreference-list.\cseq\endcsname\expandafter{\lseq, \llink{backreference.#3}{#2}}% Link: from a back reference list to the point of origin
			% Step 1: \expandafter\gdef\expandafter\"backreference-list.\cseq"\expandafter{\lseq, \llink{backreference.#3}{#2}}\fi
			% Step 2: \gdef\"backreference-list.\cseq"{\"expanded lseq", \llink{backreference.#3}{#2}}\fi
		\fi
	\fi
}
\newtoks\labelinitlist
\def\xxstripcomma, {}
\def\xstripcomma{\if\ntok,\let\xcont\xxstripcomma\else\let\xcont\relax\fi\xcont}
\def\stripcomma{\futurelet\ntok\xstripcomma}

% Prevention of duplicate labels
\inewif\ifrecorddups \recorddupstrue
\def\checkduplicates#1#2{\edef\key{\expandafter\gobble\string#1}%
	\iftypesetting\else
	\expandafter\ifx\csname line:\key\endcsname\relax\printlabel{keydefline: relax}\ifrecorddups\expandafter\xdef\csname line:\key\endcsname{\the\inputlineno}\fi\else\edef\keydefline{\csname line:\key\endcsname}\errmessage{#2}\fi\fi}

% Labels
\let\printlabel\gobble
\newtoks\curverb % curverb stores the current block name, such as "Chapter", "Theorem", etc.
\ifx\draftlabel\undefined
\def\plabel{}
\else
\def\labeltext{} % label text for proofreading purposes
\def\plabel{\ifx\labeltext\empty\else\smash{\llap{\parbackreffont\labeltext\quad}}\gdef\labeltext{}\fi} % gdef because \plabel is used inside boxes
\fi
\def\assignlabel#1#2{% #1: name (e.g., reference, paragraph, section, chapter), #2: text (e.g., 2.1)
	\ifx\lastlabel\undefined\else
	\edef\key{\expandafter\expandafter\expandafter\gobble\expandafter\string\lastlabel}% label name without backslash
	\printlabel{label \key: id.\key\space = #1.#2, text.\key\space = #2}%
	\expandafter\xdef\csname id.\key\endcsname{#1.#2}% id for DVI hrefs and backreference lists
	\expandafter\xdef\csname text.\key\endcsname{#2}% text that is actually typeset
	\edef\tmp{\the\curverb}%
	\ifx\tmp\empty\else % if curverb is nonempty, define a "verbal" label with a prefix "v"
	\printlabel{label v\key: id.v\key\space = #1.#2, text.v\key\space = \the\curverb#2}%
	\expandafter\xdef\csname id.v\key\endcsname{#1.#2}% id for DVI hrefs and backreference lists
	\expandafter\xdef\csname text.v\key\endcsname{\the\curverb#2}% text that is actually typeset
	\fi
	\predefbackref{#1.#2}%
	\fi\let\lastlabel\undefined}
\newtoks\vlist % Verification list
\def\label#1{%
	\iftypesetting\else
	\checkduplicates#1{Label \string#1 was already defined at line \keydefline}%
	\addverunused#1\verifylabel
	\def\lastlabel{#1}%
	\fi
	\numbtrue % always number labeled paragraphs
}

% Bibliography
\newif\ifyearkey % use years as bibliographic keys?
\def\y{} % use in bibliography as \y{1967}
\newdimen\bibindent % the maximum width of a bibliographic key
\newtoks\bibt % token list for all bibliographic items
\def\tbib#1{% #1 = \Paper
	\checkduplicates#1{Bibliographic reference \string#1 already defined at line \keydefline}%
	\addverunused#1\verifybib
	\edef\key{\expandafter\gobble\string#1}% reference name without backslash
	\printlabel{reference \key: id.\key\space = reference.\key, text.\key\space = \key}%
	\expandafter\edef\csname id.\key\endcsname{reference.\key}% id for DVI hrefs and backreference lists
	\expandafter\edef\csname text.\key\endcsname{\key}% text that is actually typeset
	\predefbackref{reference.\key}%
	\ifyearkey\else\setbox0=\hbox{[\key]}\ifdim\bibindent<\wd0 \bibindent=\wd0 \fi \fi % \bibindent holds the maximum length of all reference [keys]
	\appendexpand\bibt{\noexpand\typesetbib\noexpand#1\src}%
	\ifyearkey\let\next\xxftbib\else\let\next\ftbibalpha\fi\next}
\def\ftbibalpha#1\par{\append\bibt{#1\par}}
\def\xxftbib{\futurelet\next\xftbib}
\def\xftbib{\if\next[\let\next\ftbibyear\else\let\next\ftbibnoyear\fi\next}
\def\ftbibyear[#1]{\edef\yearkey{#1}\expandafter\ftbibyearbis\ignorespaces}
\def\ftbibyearbis#1\par{\append\bibt{#1\par}\ftbibend}
\def\ftbibnoyear#1\par{\append\bibt{#1\par}\edef\yearkey{\extractyear#1\par}\ftbibend}
\def\ftbibend{\expandafter\ifx\csname year.\yearkey\endcsname\relax
		\edef\yearindex{0}%
	\else
		\edef\yearindex{\csname year.\yearkey\endcsname}%
	\fi
	{\count0=\yearindex
	\advance\count0 by 1 %
	\expandafter\xdef\csname year.\yearkey\endcsname{\the\count0 }%
	\xdef\alphakey{\ifcase \count0 ?\or a\or b\or c\or d\or e\or f\or g\or h\or i\or j\or k\or l\or m\or n\or o\or p\or q\or r\or s\or t\or u\or v\or w\or x\or y\or z\else .\the\count0 \fi}}%
	\printlabel{key \key, year key \yearkey, alpha key \alphakey}%
	\expandafter\edef\csname text.\key\endcsname{\noexpand\typesetyearalpha{\key}{\yearkey}{\alphakey}}%
	\printlabel{reference \key\space adjustment: id.\key\space = reference.\key, text.\key\space = \yearkey.\alphakey}%
	\setbox0=\hbox{[\yearkey.\alphakey]}\ifdim\bibindent<\wd0 \bibindent=\wd0 \fi% \bibindent holds the maximum length of all reference [keys]
}
\def\typesetyearalpha#1#2#3{%
	\edef\yearindex{\csname year.#2\endcsname}%
	\ifnum\yearindex=1 #2\else#2.#3\fi}
\def\extractyear#1\y#2#3\par{#2}
\newif\iftype
\def\typesetbib#1#2\par{\edef\key{\expandafter\gobble\string#1}%
	\edef\bibbr{\csname backreference-list.reference.\key\endcsname}%
	\typetrue\ifsuppressunusedbib\ifx\bibbr\empty\typefalse%\warning{Suppressing unused bibliography item \string#1}
	\fi\fi
	\iftype
	\noindent\hbox to \bibindent{[\anchor{reference.\key}{\csname text.\key\endcsname}]\hfil}#2% Anchor: reference
	\ifx\bibbr\empty\else\expandafter\stripcomma\bibbr.\fi
	\hangindent\bibindent\filbreak
	\fi}

% Diagnostic messages for unused references and labels
\let\printverify\gobble
\def\addverunused#1#2{\appendexpand\vlist{\noexpand#2\noexpand#1{\the\inputlineno}}}
\def\verifyref#1#2#3{\printverify{verifying for #3 \string#1 (line #2)}%
	\edef\key{\expandafter\gobble\string#1}%
	\edef\cseq{\csname id.\key\endcsname}%
	\edef\tmp{\csname backreference-list.\cseq\endcsname}%
        \ifx\tmp\empty\ewarningline{#2}{#3 \string#1.}{\string#1}\fi}
\def\verifylabel#1#2{\verifyref#1{#2}{Unused label}}
\def\verifybib#1#2{\verifyref#1{#2}{Unused reference}}

% URLs
\let\printurl\gobble
\newtoks\urltext
\newtoks\urlt
\newif\ifpunct
\def\urldash{-}
\def\urltilde{{\tensy^^X}} % like \sim
\def\ndash{\def\urldash{--}}
\def\http://{\hfil\penalty900\hfilneg\urltext={http://}\urlt={http:/\negthinspace/}\punctfalse\urlgrab}
\def\https://{\hfil\penalty900\hfilneg\urltext={https://}\urlt={https:/\negthinspace/}\punctfalse\urlgrab}
\def\urlgrab{\catcode`\#=11 \catcode`\&=11 \futurelet\ntok\urldispatch}
\def\urldispatch{%
	\ifx\ntok~\let\proceed\urlcont\else
	\ifcat\noexpand\ntok\space\let\proceed\urlfinish\else
	\ifcat\noexpand\ntok\relax\let\proceed\urlfinish\else
	\let\proceed\urlcont
	\fi\fi\fi\proceed}
\def\urlcont#1{\ifpunct\appendexpand\urltext\punctc\appendexpand\urlt\punctc\punctfalse\fi
	\ifx\ntok~\appendexpand\urltext{\noexpand~}\appendexpand\urlt\urltilde
	\else\if\ntok\ampersand\appendexpand\urltext{&}\appendexpand\urlt{\&}%
	\else\if\ntok\ohash\appendexpand\urltext\ohash\appendexpand\urlt\#%
	\else\if\ntok_\appendexpand\urltext_\appendexpand\urlt\_%
	\else\if\ntok-\appendexpand\urltext-\appendexpand\urlt\urldash
	\else\if\ntok.\puncttrue\def\punctc{.}%
	\else\if\ntok,\puncttrue\def\punctc{,}%
	\else\if\ntok;\puncttrue\def\punctc{;}%
	\else\appendexpand\urltext{#1}\appendexpand\urlt{#1}%
	\fi\fi\fi\fi\fi\fi\fi\fi\urlgrab}
\def\urlfinish{\catcode`\#=6 \catcode`\&=4 \hbox{\printurl{\the\urltext}\link{\the\urltext}{\the\urlt}}\ifpunct\punctc\punctfalse\fi\def\urldash{-}}
\def\idgrab{\futurelet\ntok\iddispatch}
\def\iddispatch{\ifcat\noexpand\ntok\space\let\proceed\urlfinish
		\else\if\ntok,\let\proceed\urlfinish
		\else\let\proceed\idcont
		\fi\fi\proceed}
\def\idcont#1{\ifpunct\appendexpand\urltext.\appendexpand\urlt.\punctfalse\fi
	\if\ntok.\puncttrue\def\punctc{.}%
	\else\if\ntok_\appendexpand\urltext_\appendexpand\urlt\_%
	\else\appendexpand\urltext{#1}\appendexpand\urlt{#1}%
	\fi\fi\idgrab}

% Mathematical string grabbing
\let\printgrab\gobble
\newtoks\grabname
\newtoks\grabtoks % beginning
\newtoks\grabcseq % control sequence name
\newtoks\subsuptoks % ^ or _
\newtoks\dtoks % first letter of a subscript
\newcount\grabsize
\newif\ifgrabsubscript % include subscript?
\newif\ifgrabsupscript % include superscript?
\def\grabsequence{\bgroup % \bgroup allows us to say things like $C^@op$, with "op" being a superscript
	\grabsubscripttrue\grabsupscripttrue % do grab sub/superscripts
	\grabstring}
\def\grabalpha{\bgroup % same
	\grabsubscriptfalse\grabsupscriptfalse % do not grab sub/superscripts
	\grabstring}
\def\grabingroup{\ifinfont\errmessage{Already inside a math token}\fi\append\grabtoks{\bgroup\grablink}\infonttrue}
\def\graboutgroup{\ifinfont\append\grabtoks{\endgrablink\egroup}\infontfalse\fi}
\def\grabstring#1#2#3{% #1 = descriptive name like cat, fun, trans; #2 = font like \bf, \rm, \it; #3 = postcommand like \nolimits
	\let\specialhat^
	\catcode`\^=7
	\aftergroup#3 % #3 could be \nolimits or \limits
	\ifx\specialaddon\undefined\else\expandafter\aftergroup\specialaddon\let\specialaddon\undefined\fi
	\inewif\ifdefine \inewif\ifinfont
	\printgrab{}\printgrab{grab a string of type #1, typeset using font \string#2, with postcommand \string#3}%
	\grabname{#1}\def\grabfont{#2}\grabsize0 \grabtoks={}\grabingroup \append\grabtoks{#2}\grabcseq={}%
	\futurelet\ntok\grabdeflookahead}
\def\grabdeflookahead{\if=\noexpand\ntok % @=Set creates an anchor, whereas @Set refers to it
	\definetrue\printgrab{defining}\expandafter\grabgobblefuturelet
	\else\printgrab{referencing}\definefalse\expandafter\grablookahead\fi}
\def\grabgobblefuturelet#1{\futurelet\ntok\grabtestforsilent} % gobble = and look for another =
\newif\ifsilentgrab
\def\grabtestforsilent{\if=\noexpand\ntok \silentgrabtrue \let\ncom\grabsilenteq \else \silentgrabfalse \let\ncom\grablookahead \fi \ncom}
\def\grabsilenteq={\grabfuturelet}
\def\grabfuturelet{\futurelet\ntok\grablookahead}
\def\grablookahead{\printgrab{futurelet token meaning: \meaning\ntok}%
	\let\ncom\grabfinish
	\if\bgroup\noexpand\ntok \printgrab{left brace, terminating}%
	\else \if\egroup\noexpand\ntok \printgrab{right brace, terminating}%
	\else \if\space\noexpand\ntok \printgrab{blank space, terminating}%
	\else \let\ncom\grabexamine \fi\fi\fi \ncom}
\def\grabexamine#1{\printgrab{grabexamine argument: \string#1, meaning \meaning#1}%
	\def\ncom{\grabfinish#1}%
	\ifcat$\ifcat*\string#1\fi$% is #1 not a command sequence?
		\ifcat _\noexpand#1 \ifgrabsubscript\printgrab{subscript, continuing}%
			\graboutgroup \append\grabtoks{#1}\subsuptoks{#1}\def\ncom{\grabsubsupfuturelet}%
						\else\printgrab{subscript, terminating}\fi
		\else \ifcat ^\noexpand#1 \ifgrabsupscript\printgrab{superscript, continuing}%
			\graboutgroup \append\grabtoks{#1}\subsuptoks{#1}\def\ncom{\grabsubsupfuturelet}%
						\else\printgrab{superscript, terminating}\fi
		\else \ifx \specialhat#1 \ifgrabsupscript\printgrab{specialhat superscript, continuing}%
			\graboutgroup \append\grabtoks{#1}\subsuptoks{#1}\def\ncom{\grabsubsupfuturelet}%
						\else\printgrab{superscript, terminating}\fi
		\else \ifcat\noexpand~\noexpand#1 \printgrab{active character \string#1, examining further}%
			\ifnum1=\uccode`#1 \printgrab{UTF-8 letter, continuing}%
				\advance\grabsize1 \append\grabtoks{#1}\appendexpand\grabcseq{\string#1}\def\ncom{\grabfuturelet}%
			\else
				\ifnum\the\grabsize=0 \printgrab{Nothing grabbed so far, continuing}%
					\advance\grabsize1 \append\grabtoks{#1}\appendexpand\grabcseq{\string#1}\def\ncom{\grabfuturelet}%
				\else\printgrab{Not a UTF-8 letter and not the first character in a string, terminating}%
				\fi
			\fi
		\else \ifcat a\noexpand#1 \printgrab{letter #1, continuing}%
			\advance\grabsize1 \append\grabtoks{#1}\append\grabcseq{#1}\def\ncom{\grabfuturelet}%
		\else\printgrab{nonactive character \string#1}%
			\ifnum\the\grabsize=0 \printgrab{sole argument, adding and terminating}%
				\advance\grabsize1 \append\grabtoks{#1}\append\grabcseq{#1}\def\ncom{\grabfinish}%
			\else\printgrab{terminating}\fi
		\fi\fi\fi\fi\fi
	\else \printgrab{command sequence \string#1, terminating}\fi
	\ncom}
\def\grabsubsupfuturelet{\futurelet\ntok\grabsubsuplookahead}
\newcount\dig
\newif\ifdigit
\def\grabsubsuplookahead{\printgrab{subsup futurelet token meaning: \meaning\ntok}%
	\if\bgroup\noexpand\ntok \printgrab{left brace, continuing}\let\ncom\grabentiresubsup%
	\else \if\egroup\noexpand\ntok \printgrab{right brace, continuing}\let\ncom\grabentiresubsup%
	\else \if\space\noexpand\ntok \errmessage{Blank space after \the\subsuptoks}%
	\else \let\ncom\grabsubsupexamine \fi\fi\fi \ncom}
\def\grabentiresubsup#1{\printgrab{subsup entire group added}\grabingroup\append\grabtoks{#1}\graboutgroup\grabfuturelet}
\def\grabsubsupexamine#1{\printgrab{examining subsup argument \string#1, meaning \meaning#1}%
	% We pass through (1) single letters; (2) command sequences
	\ifcat$\ifcat*\string#1\fi$% is #1 not a command sequence?
		\ifcat\noexpand~\noexpand#1 \printgrab{active character \string#1, continuing}%
			%\warning{math active C: \string#1, grabsize=\the\grabsize, grabtoks=\the\grabtoks, grabcseq=\the\grabcseq}%
			%\grabingroup\appendexpand\grabtoks{\grabfont\noexpand#1}\let\ncom\grabsubsupremainderfuturelet
			\grabingroup\appendexpand\grabtoks{\grabfont\noexpand#1}\let\ncom\grabfuturelet
		\else\ifnum"8000=\the\mathcode`#1 \printgrab{math active character \string#1, continuing}%
			%\warning{math active A: \string#1, grabsize=\the\grabsize, grabtoks=\the\grabtoks, grabcseq=\the\grabcseq}%
			%\grabingroup\appendexpand\grabtoks{\grabfont\noexpand#1}\let\ncom\grabsubsupremainderfuturelet
			\let\specialaddon\egroup
			\def\ncom{#1}%
			\grabtypeset
		\else\ifcat a\noexpand#1 \printgrab{letter #1, checking whether single or not}%
			\dtoks{#1}\def\ncom{\futurelet\ntok\grabsubsupsecondletterlookahead}%
		\else\printgrab{something else, inserting a single-character sub/superscript, continuing}
			\appendexpand\grabtoks{\bgroup\grabfont\noexpand#1\egroup}%
			\advance\grabsize2 \appendexpand\grabcseq{\the\subsuptoks\string#1}% possibly ignore digits here
			\def\ncom{\grabfuturelet}\fi\fi\fi
	\else \printgrab{command sequence \string#1, continuing}%
		\append\grabtoks{#1}\let\ncom\grabfuturelet\fi
	\ncom}
\def\grabsubsupsecondletterlookahead{\def\ncom{\appendexpand\grabtoks{\the\dtoks}\grabfuturelet}%
	\ifcat a\noexpand\ntok \printgrab{not a single letter, grabbing the entire subsupscript}%
		\advance\grabsize1 \appendexpand\grabcseq{\expandafter\string\the\subsuptoks}% append _ or ^ to the label
		\grabingroup
		\appendexpand\grabtoks{\grabfont\the\dtoks}%
		\appendexpand\grabcseq{\the\dtoks}%
		\def\ncom{\grabfuturelet}%
	\else \printgrab{single letter, continuing}\fi\ncom}
\def\grabfinish{\printgrab{grabfinish}\graboutgroup\grabtypeset\egroup}
\def\grabtypeset{\printgrab{grabtypeset grabsize=\the\grabsize, grabtoks=\the\grabtoks, grabcseq=\the\grabcseq}%
	\def\grablink##1\endgrablink{##1}%
	\ifnum\the\grabsize=0 \errmessage{No string to grab}\fi
	\ifnum\the\grabsize>1 % single-letter names are not references
		\ifdefine % defining a mathematical identifier
			\expandafter\checkduplicates\csname\the\grabname.\the\grabcseq\endcsname{Mathematical identifier \key\space already defined at line \keydefline}%
			\iftypesetting % if we are actually typesetting, create an anchor
				\ifsilentgrab
					\expandafter\gdef\csname silent:\the\grabname.\the\grabcseq\endcsname{}% record that this id is silent
				\else
					\anchor{\the\grabname.\the\grabcseq}{}% Anchor: definition of a mathematical identifier
				\fi % only create if not silent
			\else
				\expandafter\xdef\csname id.\the\grabname.\the\grabcseq\endcsname{paragraph.\the\secn.\the\parn}% paragraph id for a back reference list
				\predefbackref{paragraph.\the\secn.\the\parn}%
  			\fi
		\else % referencing a mathematical identifier
			\iftypesetting % we are actually typesetting
				\global\advance\backref1
				\expandafter\ifx\csname line:\the\grabname.\the\grabcseq\endcsname\relax % identifier is undefined
					\warning{Undefined mathematical identifier \the\grabname.\the\grabcseq}%
					\expandafter\gdef\csname\the\grabname.\the\grabcseq\endcsname{\relax}% report undefined references only once; gdef because inside bgroup..egroup
				\else % identifier has ben defined
					\expandafter\ifx\csname silent:\the\grabname.\the\grabcseq\endcsname\empty
						\edef\grablink##1\endgrablink{{##1}}%
					\else % not silent, need a hyperlink
						\edef\grablink##1\endgrablink{\noexpand\llink{\the\grabname.\the\grabcseq}{##1}% Link: from a mathematical identifier to its definition 
							\noexpand\anchor{backreference.\the\backref}{}}% Anchor: back reference for a mathematical identifier reference
					\fi
				\fi
			\else % not yet typesetting, just collecting back references
				\ifsuppressbackref\else
					\global\advance\backref1
					\printbackref{math back reference to \the\grabname.\the\grabcseq: backref.\the\backref\space at line \the\inputlineno}%
					\expandafter\recordbackref\csname\the\grabname.\the\grabcseq\endcsname
				\fi
			\fi
		\fi
	\fi
	\ifsilentgrab\else\the\expandafter\grabtoks\fi}
%		\ifcat 0\noexpand\ntok \printgrab{subsupscript is followed by a catcode 12 character, examining further}%
%			\dig=0 \loop \if\the\dig\ntok \digittrue \fi \ifnum\dig<9 \advance\dig1 \repeat
%			\ifdigit \printgrab{subsupscript is followed by a digit, passed through}%
%				\advance\grabsize2
%				\appendexpand\grabtoks{\subsup#1}%
%				\appendexpand\grabcseq{\expandafter\string\subsup#1}%
%				\def\ncom{\grabfinish}%
%			\else \printgrab{subsupscript is followed not by a digit or a letter, terminating}\fi

% index macros modeled after manmac.tex
% ^={...}: setup an anchor and typeset in italic font
% ^=:{...}: setup an anchor and typeset in roman font
% ^^={...}: setup an anchor, do not typeset
% ^{...}: record a back reference and typeset
% ^^{...}: record a back reference, do not typeset
% ^!{...}: reference a proclaimed statement in the format Theorem~3.5
\inewif\ifsilent \inewif\ifanchor
\newif\ifsuppresscs
\catcode`\^=7   \def\specialhat{\ifmmode\def\next{^}\else\let\next\beginxref\fi\next} \catcode`\^=\active \let^=\specialhat
\def\silentxref#1{\futurelet\next\silentxrefswitch}
\def\silentxrefswitch{\silenttrue\xref}
\def\beginxref{\futurelet\next\beginxrefswitch}
\def\beginxrefswitch{\ifx\next\specialhat\let\next\silentxref \else\silentfalse\let\next\xref\fi \next}
\def\xref{\leavevmode\futurelet\next\xrefswitch}
\def\xrefswitch{\ifx\next!\let\next\verbalxref \else \ifx\next=\let\next\anchorxref \else \anchorfalse \let\next\normalxref \fi \fi \next}
\newtoks\vtoksl
{\count0="C2 \loop\ifnum\count0<"F5 \catcode\count0=11 \advance\count0 by 1 \repeat % make all legitimate initial UTF-8 octets letters
\gdef\plainaccents{\suppresscstrue%
	\def\`##1{##1\empty ̀}%
	\def\'##1{##1\empty ́}%
	\def\^##1{##1\empty ̂}%
	\def\"##1{##1\empty ̈}%
	\def\~##1{##1\empty ̃}%
	\def\=##1{##1\empty ̄}%
	\def\.##1{##1\empty ̇}%
	\def\u##1{##1\empty ̆}%
	\def\v##1{##1\empty ̌}%
	\def\H##1{##1\empty ̋}%
	\def\t##1{##1\empty ͡}%
}}
\def\plainaccents{\let\xcsname=\empty \let\xendcsname=\empty}
\def\verbalxref!{\begingroup\plainaccents\verbalxrefaux}
\def\verbalxrefaux#1{%
	\lowercase{\vtoksl{#1}}%
	\expandtoks\vtoksl
	\iftypesetting
		\expandafter\gdef\expandafter\cseq\expandafter{\csname verbal.\the\vtoksl\endcsname}%
		\printlabel{verbal xref \the\vtoksl}%
		\expandafter\ifx\cseq\relax\errmessage{Undefined reference to \the\vtoksl}\else\cseq\fi
	\else
		\ifsuppressbackref\else
			\global\advance\backref1 %
			\printbackref{label verbal.#1: backref.\the\backref\space at line \the\inputlineno}%
			\blah
			\expandafter\recordbackref\csname verbal.\the\vtoksl\endcsname
		\fi
	\fi
	\endgroup}
\def\initverballabelcommand#1{%
	\printlabel{initializing verbal label #1 (\the\curverb\the\secn.\the\parn)}%
	\expandafter\xdef\csname id.verbal.#1\endcsname{paragraph.\the\secn.\the\parn}%
	\expandafter\xdef\csname text.verbal.#1\endcsname{\the\curverb\the\secn.\the\parn}%
	\expandafter\initlabelcommand\csname verbal.#1\endcsname
}
\def\anchorxref={\anchortrue\futurelet\next\anchorxrefswitch}
\def\anchorxrefswitch{\ifx\next:\let\next\nonitalicanchor\else\italictrue\let\next\normalxref\fi \next}
\def\nonitalicanchor:{\italicfalse\normalxref}
\newtoks\firsttoks \newtoks\secondtoks \inewif\ifplural \inewif\ifitalic
\def\parseplural#1[#2|#3]{\let\next\parseplural\ifx\hfuzz#2\hfuzz\ifx\hfuzz#3\hfuzz\let\next\relax\else\pluraltrue\fi\else\pluraltrue\fi
	\append\firsttoks{#1#2}\append\secondtoks{#1#3}\next}
\newtoks\firsttoksl
\newtoks\secondtoksl
\newtoks\nexttoks
\def\normalxref{\begingroup\plainaccents\normalxrefaux}
\def\normalxrefaux#1{\firsttoks{}\secondtoks{}\pluralfalse\parseplural#1[|]%
	\lowercase\expandafter{\expandafter\firsttoksl\expandafter{\the\firsttoks}}%
	\expandtoks\firsttoksl
	\lowercase\expandafter{\expandafter\secondtoksl\expandafter{\the\secondtoks}}%
	\expandtoks\secondtoksl
	\ifanchor
		\iftypesetting\else
			\initverballabelcommand{\the\firsttoksl}%
			\ifplural\initverballabelcommand{\the\secondtoksl}\fi
			\predefbackref{paragraph.\the\secn.\the\parn}%
			\expandafter\checkduplicates\csname\the\firsttoksl\endcsname{Verbal label \the\firsttoksl\space was already defined at line \keydefline}
			\ifplural\expandafter\checkduplicates\csname\the\secondtoksl\endcsname{Verbal label \the\secondtoksl\space was already defined at line \keydefline}\fi
		\fi
		\anchor{verbal.\the\firsttoksl}{}% Anchor: verbal reference
		\ifplural\anchor{verbal.\the\secondtoksl}{}\fi % Anchor: plural verbal reference
	\else
		\ifsuppressbackref\else
			\global\advance\backref1
			\printbackref{label normal.#1: backref.\the\backref\space at line \the\inputlineno}%
			\iftypesetting
			\else
				\expandafter\recordbackref\csname verbal.\the\secondtoksl\endcsname
			\fi
			\anchor{backreference.\the\backref}{}% Anchor: back reference for a verbal reference
		\fi
	\fi
	% The actual text is typeset after \endgroup
	\nexttoks{}%
	\ifsilent
		\nexttoks{\ignorespaces}%
	\else
		\ifanchor
			\ifitalic\nexttoks\expandafter{\expandafter\bgroup\expandafter\it\the\firsttoks\italcorr}%
			\else\nexttoks\expandafter{\the\firsttoks}%
			\fi
		\else
			\edef\tmp{{verbal.\the\secondtoksl}}%
			\nexttoks\expandafter{\expandafter\llink\tmp}% Link: from a verbal reference to its definition
			\nexttoks\expandafter\expandafter\expandafter{\expandafter\the\expandafter\nexttoks\expandafter{\the\firsttoks}}%
		\fi
	\fi
	\expandafter\endgroup\the\nexttoks}
\def\italcorr{\futurelet\next\italcorrtest}
\def\italcorrtest{\if,\noexpand\next\else\if.\noexpand\next\else\/\fi\fi\egroup}

% Individual repositories

%\def\mathjournals{\ndash\http://www.mathjournals.org}
%\def\numdam{\urltext={http://www.numdam.org/item/?id=}\urlt={numdam:}\punctfalse\expandafter\idgrab\ignorespaces}

%\def\matwbn{\http://matwbn.icm.edu.pl/ksiazki}
\def\gen:{\http://libgen.rs/book/index.php?md5=}
\def\jstor:{\https://www.jstor.org/stable/}
\def\eudml:{\https://eudml.org/doc/}

\def\arXiv:{\urltext={https://arxiv.org/abs/}\urlt={arXiv:}\punctfalse\idgrab}

\def\Zbl:{\urltext={https://zbmath.org/?q=an:}\urlt={Zbl:}\punctfalse\idgrab}
\def\doi:{\ndash\urltext={https://doi.org/}\urlt={doi:}\punctfalse\urlgrab}

% Squares
\def\sqr#1#2{{\thinspace\vbox{\hrule height.#2pt \hbox{\vrule width.#2pt height#1pt \kern#1pt \vrule width.#2pt} \hrule height0pt depth.#2pt}\thinspace}}
\def\square{\mathchoice\sqr64\sqr64\sqr{4.2}3\sqr33}

% Long arrows
\def\ltoarr#1{\mathop{\count0=#1 \loop\ifnum\count0>0 \smash-\mkern-7mu \advance\count0 -1 \repeat \mathord\rightarrow}\limits} % parametrized \rightarrowfill
\def\lto#1#2{\mathrel{\ltoarr{#1}^{#2}}} % parametrized \rightarrowfill, with a label
\def\longto#1^#2_#3{\mathrel{\ltoarr{#1}^{#2}_{#3}}} % parametrized \rightarrowfill, with a label above and below
\def\lgetsarr#1{\mathop{\mathord\leftarrow \count0=#1 \loop\ifnum\count0>0 \mkern-7mu\smash-\advance\count0 -1 \repeat}\limits} % parametrized \leftarrowfill
 % parametrized \leftarrowfill, with a label
\def\longgets#1^#2_#3{\mathrel{\lgetsarr{#1}\limits^{#2}_{#3}}} % parametrized \leftarrowfill, with a label

% Double arrows

\def\toto{\mathrel{\vcenter{\hbox{$\to$}\kern-1.5ex \hbox{$\to$}}}}
\def\prearrfill{\smash-\mkern-7mu}
\def\postarrfill{\mkern-7mu\smash-}
\def\midarrfill#1{\cleaders\hbox{$\mkern-2mu\smash-\mkern-2mu$}\hskip0pt plus #1fil}
\def\rightarrfill{\mkern-7mu\mathord\rightarrow}
\def\leftarrfill{\mathord\leftarrow\mkern-7mu\midarrfill1\postarrfill}
\def\ltoto#1#2#3{\ifinner
	\mathrel{\vcenter{\hbox to #1em{$\prearrfill\midarrfill1{\scriptstyle#2}\midarrfill3 \rightarrfill$}%
		\kern-1.5ex \hbox to #1em{$\prearrfill\midarrfill3{\scriptstyle#3}\midarrfill1 \rightarrfill$}}}%
	\else
	\mathrel{\mathop{\vcenter{\hbox to #1em{\rightarrowfill}%
		\kern-1.5ex \hbox to #1em{\rightarrowfill}}}\limits^{#2}_{#3}}%
	\fi}
\def\ltogets#1#2#3{\ifinner
	\mathrel{\vcenter{\hbox to #1em{$\prearrfill\midarrfill1{\scriptstyle#2}\midarrfill3 \rightarrfill$}%
		\kern-1.5ex \hbox to #1em{$\leftarrfill\midarrfill3{\scriptstyle#3}\midarrfill1 \postarrfill$}}}%
	\else
	\mathrel{\mathop{\vcenter{\hbox to #1em{\rightarrowfill}%
		\kern-1.5ex \hbox to #1em{\leftarrowfill}}}\limits^{#2}_{#3}}%
	\fi}

% New macros for arrows
\def\ltogetscore#1#2{\dimen0=\fontdimen6 #1 2 \divide\dimen0 by 2 \multiply\dimen0 by #2 \vcenter{\hbox to \dimen0{\rightarrowfill}\kern-1.8ex \hbox to \dimen0{\leftarrowfill}}}
\def\ltogets#1#2#3{\mathrel{\mathop{\mathchoice{\ltogetscore\textfont{#1}}{\ltogetscore\textfont{#1}}{\ltogetscore\scriptfont{#1}}{\ltogetscore\scriptscriptfont{#1}}}^{#2}_{#3}}}

% Diagonal arrows
\def\rx#1#2{\rlap{\kern #1pt \raise#1pt \hbox{#2}}}
\def\dottednearrow{\rx{-8}. \rx{-6}. \rx{-4}. \rx{-2}. \rx0. \rx2. \rx4. \kern6pt \raise7.7pt \hbox{$\nearrow$}}

% Rectangular commutative diagrams
% \matrix from Plain TeX: #1 = \hfil$##$\hfil&&\quad\hfil$##$\hfil
\def\gmatrix#1#2{\null\,\vcenter{\normalbaselines
	\ialign{#1\crcr
		\mathstrut\crcr\noalign{\kern-\baselineskip}
		#2\crcr\mathstrut\crcr\noalign{\kern-\baselineskip}}}\,}

\def\sqmatrix{\gmatrix{\hfil$##$&\enspace\hfil$##$\hfil\enspace&$##$\hfil}}
\def\cdbl{\def\normalbaselines{\baselineskip20pt \lineskip3pt \lineskiplimit3pt }}

\def\sqcd{\cdbl\let\vagap\;\sqmatrix}

% Horizontal and vertical diagram arrows
\newcount\arrowsize \arrowsize3
\def\mapright#1{\smash{\lto\arrowsize{#1}}}

\def\rvagap{\vagap}  \def\rvaskip{\vaskip}  \def\vaskip{} \def\vagap{}
\def\mapdown#1{\rvagap\Big\downarrow\rlap{$\vcenter{\hbox{$\scriptstyle#1$}}$}\rvaskip}

% Wrappable diagrams
\newcount\forno \forno0
\def\arrno#1#2{\global\advance\arr1 \edef\eeqnno{\the\arr}%
	\global\advance\forno1 \edef\eforno{\the\forno}%
	\xdef#2{\noexpand\llink{equation.\eforno}{\eeqnno}}% Link: from an equation number to its definition; xdef to ensure it can be seen outside of \wrapdiagram
	#1{(\anchor{equation.\eforno}{\eeqnno})}} % Anchor: equation number
\newbox\mdiag
\def\wrapdiagram{%
	\setbox\mdiag\vtop\bgroup
	\null % an empty hbox to ensure proper vertical positioning
	\vskip\baselineskip
	\inewcount\arr \arr0
	\baselineskip0pt
	\lineskip4pt
	\lineskiplimit4pt
	\let\par\cr
	\obeylines
	\halign\bgroup\hfil$\displaystyle##$\hfil\cr
	\ewrapdiagram}

\def\ewrapdiagram#1{#1
	\egroup
	\egroup
	\vskip0pt plus \dp\mdiag \penalty-250 \vskip0pt plus-\dp\mdiag % ensure the diagram fits on a single page
	\hangafter-\dp\mdiag
	\divide\hangafter\baselineskip
	\advance\hangafter-2
	\hangindent-\wd\mdiag
	\advance\hangindent-2em
	\hbox to\hsize{\hfil\dp\mdiag0pt \box\mdiag}%
	\ignorespaces}

% METAPOST diagrams                                                                                                     
\newlinechar=10
\inewtoks\preamble
{\endlinechar=10 \catcode`#=12 \global\preamble{
prologues := 3;

verbatimtex
\let\endprolog

\expandafter\gobbleinit\input }\global\appendonceexpand\preamble\jobname\global\append\preamble{
\catcode`\^=7
etex

input cmarrows
setup_cmarrows(arrow_name = "texarrow"; parameter_file = "cmr10.mf"; macro_name = "drawarrow");
setup_cmarrows(arrow_name = "doublearrow"; parameter_file = "cmr10.mf"; macro_name = "drawdarrow");
def drawmarrow expr p = _apth:=p; _finmarr enddef;
rule_thickness#:=.4pt#;    % cmr10.mf: thickness of lines in math symbols
def _finmarr text t_ =
  drawarrow subpath(0, 0.5 * length(_apth)) of _apth t_;
  draw subpath(0.5 * length(_apth), length(_apth)) of _apth withpen pencircle scaled rule_thickness# t_;
enddef;

def object(suffix O)(expr x,y)(expr l) =
  save O;       
  pair O;
  O := (x,y) * u;
  picture O.tx;    
  O.tx := thelabel(l,O);
  draw O.tx;                         
enddef;                                    

def smorphism(suffix A,B) =
  save ss, tt;
  ss := xpart ((A..B) intersectiontimes bbox A.tx);
  tt := xpart ((A..B) intersectiontimes bbox B.tx);
  drawarrow subpath(ss,tt) of (A..B);
enddef;

def morphism(suffix A,B)(expr l)(expr f)(suffix $) =
  smorphism(A, B);
  label.$(l, point f[ss, tt] of (A..B));
enddef;
}}%

\newcount\vertex                                      
\ifx\mathspecials\undefined\def\mathspecials{}\fi

\def\grabdiagramaux#1;{\vertex0 \dcom#1,*}                          
\def\dcom{\futurelet\next\dcomswitch}                                          
\def\dcomswitch{\ifcat a\noexpand\next \let\next\dcomalpha                                       
        \else\if*\noexpand\next \addcode{\the\buffertoks}\endinlinemp\let\next\gobble                                                 
        \else\if[\noexpand\next \let\next\grabscale                                                                
        \else\errmessage{Unrecognized diagram command \next}\let\next\relax\fi\fi\fi\next}
\def\grabscale[#1]{\dimen0=#1 \addcode{save u; u = \the\dimen0;}\dcom}
\def\dcomalpha#1{\def\objectname{#1}\futurelet\next\dcomalphaswitch}
\def\dcomalphaswitch{\ifcat a\noexpand\next \let\next\dcommorphism\else                                                          
        \expandafter\edef\csname vertex:\objectname\endcsname{\the\vertex}%              
        \if:\noexpand\next\let\next\grabcoords                                              
        \else\if=\noexpand\next\let\next\grabobjectequ                                     
        \else\errmessage{Expected : or = while processing a diagram object, got \meaning\next}\let\next\relax\fi\fi\fi\next}
\def\grabcoords:#1,#2={\toks\vertex{#1,#2}\grabobjectlabel}
\def\grabobjectequ={\edef\tmp{\the\toks\vertex}\ifx\tmp\empty\errmessage{No coordinates specified for vertex \the\vertex: \objectname}\fi\grabobjectlabel}
\def\grabobjectlabel#1,{\addcode{object(\objectname, \the\toks\vertex, }\addunexpandedcode{btex #1 etex);}\advance\vertex1 \dcom}
\def\dcommorphism#1{\def\tobjectname{#1}%
        \def\labelpos{.5}%                                      
        \futurelet\next\dcommorphismswitch}                     
\def\dcommorphismswitch{\if.\noexpand\next \expandafter\grabmorphismdir \else \setupmorphismlabel \expandafter\grabmorphismpos\fi}
\def\setupmorphismlabel{
        \edef\vlabeldira{direction.\objectname.\tobjectname}%
        \edef\vlabeldirb{direction.\tobjectname.\objectname}%
        \expandafter\ifx\csname\vlabeldira\endcsname\relax
        	\expandafter\ifx\csname\vlabeldirb\endcsname\relax
		        \edef\tmpa{\csname vertex:\objectname\endcsname}\expandafter\ifx\tmpa\relax\errmessage{No such vertex: \objectname}\fi
		        \edef\tmpb{\csname vertex:\tobjectname\endcsname}\expandafter\ifx\tmpb\relax\errmessage{No such vertex: \tobjectname}\fi
		        \ifnum \tmpa<\tmpb \edef\vlabeldir{direction.\tmpa.\tmpb}\else \edef\vlabeldir{direction.\tmpb.\tmpa}\fi
		        \expandafter\ifx\csname\vlabeldir\endcsname\relax\errmessage{No label direction specified for morphism \objectname->\tobjectname}\fi
		        \edef\labeldir{\csname\vlabeldir\endcsname}%                        
		\else
        		\edef\labeldir{\csname\vlabeldirb\endcsname}%                        
		\fi
	\else
        	\edef\labeldir{\csname\vlabeldira\endcsname}%                        
	\fi
}                                                                                        
\def\grabmorphismdir.{\def\labeldir{}\futurelet\next\grabmorphismdirec}
{\catcode`\@=13                                                               
\gdef\grabmorphismdirec{\ifx@\next \let\next\grabmorphismposition
        \else\if=\noexpand\next \let\next\grabmorphismequ               
        \else \let\next\grabmorphismdirect\fi\fi\next}                       
\gdef\grabmorphismdirect#1{\edef\labeldir{\labeldir#1}\futurelet\next\grabmorphismdirec}
\gdef\grabmorphismpos{\ifx@\next \expandafter\grabmorphismposition \else \expandafter\grabmorphismequ\fi} 
\gdef\grabmorphismposition@#1={\def\labelpos{#1}\grabmorphismlabel}
}
\def\grabmorphismequ={\grabmorphismlabel}
\def\grabmorphismlabel#1,{%
        \addcodebuffer{morphism(\objectname, \tobjectname, }%
        \addunexpandedcodebuffer{btex \everymath{\scriptstyle}#1 etex, }%
        \addcodebuffer{\labelpos, \labeldir);}%
        \dcom}

% Label extraction and definition
\let\printextract\gobble
\def\hinitlabelcommand#1{\printextract{initializing label command \string#1}%
	\gdef#1{\printlabel{invoked label \string#1}% #1: label command, #2: label id for DVI, #3: typesetted text
		\ifsuppressbackref
			\edef\key{\expandafter\gobble\string#1}%
			\csname text.\key\endcsname
		\else
			\global\advance\backref1 %
			\printbackref{label \expandafter\gobble\string#1: backref.\the\backref\space at line \the\inputlineno}%
			\iftypesetting
			\else
				\blah
				\recordbackref#1
			\fi
			\anchor{backreference.\the\backref}{}% Anchor: back reference for a label (theorem or bibliography)
			\edef\key{\expandafter\gobble\string#1}%
			\llink{\csname id.\key\endcsname}{\csname text.\key\endcsname}% Link: from a label (theorem or bibliography) to its definition
		\fi
		}}
\let\initlabelcommand\hinitlabelcommand
\def\pinitlabelcommand#1{\printextract{initializing label command \string#1}%
	\gdef#1{\printlabel{invoked plain label \string#1}% #1: label command, #2: label id for DVI, #3: typesetted text
		\iftypesetting
			\edef\key{\expandafter\gobble\string#1}%
			\csname text.\key\endcsname
		\else
			\blah
		\fi}}
\newread\labelin
\newif\iflabelcont
\let\terminate=\relax % allow \terminate to be read by \read
\long\def\labelauxaux#1\terminate
{}
\def\preprocesslabel#1{%
	\printextract{Processing label \string#1}%
	\edef\key{\expandafter\gobble\string#1}%
	\initlabelcommand#1
	\expandafter\initlabelcommand\csname v\key\endcsname
	\labelauxaux
}
\def\preprocessbib#1{%
	\printextract{Processing bib \string#1}%
	\initlabelcommand#1
	\labelauxaux
}
\long\def\labelaux#1{\ifx#1\label\let\next\preprocesslabel\else\ifx#1\bib\let\next\preprocessbib\else\let\next\labelauxaux\fi\fi\next}
\def\processoneline{\expandafter\labelaux\labelline\relax\relax\relax\terminate}
\def\preprocess#1{% preprocessing stage to define all labels, back references, bibliographic items, and the table of contents
	% first pass: predefine command sequences for labels and references
	\openin\labelin=#1 \labelconttrue
	\loop
		\read\labelin to\labelline
		\ifeof\labelin\let\next\labelcontfalse\else\let\next\processoneline\fi
		\next
	\iflabelcont\repeat
	% second pass: process labels, bibliography, back references, and the table of contents
	\expandafter\gobbleinit\input#1
	\ifhmode\par\fi\vfill\eject % for LaTeX
	\ifhmode\par\fi\vfill\supereject
}

\def\importlabels#1{ % import labels from an external document
	\recorddupsfalse
	\let\initlabelcommand\pinitlabelcommand
	\preprocess{#1}%
	\let\initlabelcommand\hinitlabelcommand
	\recorddupstrue
	\cont={}% discard the external table of contents
	\vlist={}% discard the verification list
	\let\chapname\undefined \secn0 \parn0 \backref0 % reset the numbering
}

% Primary typesetting routine

\newif\iftypesetting % are we in the final typesetting stage?

% Label and bibliography collection stage
\typesettingfalse % no typesetting at this stage
\def\blah{blah} % placeholder for future references
\output{\setbox0\box255 \setbox0\box\footins \deadcycles0 } % fake typesetting
\let\bib\tbib % collect bibliographic items at this stage
\let\refs\relax % do not typeset bibliography
\everypar{\numpar}
\parn0
\def\par{\finishpar} % no paragraph back references
% Collect embedded METAPOST files and write them to \jobname.gen.mp
\newif\ifmetapost \metapostfalse
\edef\filestem{\jobname.gen}
\newcount\figno
\newwrite\mpout
\newtoks\outtoks
\def\inimp{\ifmetapost\else\global\metaposttrue\immediate\openout\mpout=\filestem.mp \addcode{\the\preamble}\fi}
\figno0
% Temporarily set up large width and height to suppress over/underfull boxes
\newdimen\oldhsize
\oldhsize=\hsize
\newdimen\oldvsize
\oldvsize=\vsize
\hsize\maxdimen
\vsize\maxdimen
\hbadness10000
% Preprocess the file
\preprocess\jobname % collect all labels and bibliography from the main documents and possibly also external labels
% Restore the previously saved dimensions
\hsize\oldhsize
\vsize\oldvsize
\ifx\fmtname\plainfmtname
\hsize\plaintextwidth
\vsize\plaintextheight
\fi
\hbadness1000
% Process back references
\printerr{(\jobname.tex}
\printbackref{list of all back references: \the\backreflist}%
\the\backreflist % process all back references
\ifsuppressunusedbib\def\verifybib#1#2{}\fi
\the\vlist % verify if there are any unused references or bibliography items
\printerr{)}
% Finalize and compile the METAPOST file
\ifmetapost
\addcode{end}
\immediate\closeout\mpout
\immediate\write18{mpost -interaction nonstopmode \filestem.mp}% can also compile the METAPOST file separately
\ifeof18 \warning{Compile the METAPOST file \filestem.mp manually using mpost \filestem.mp}\let\runmp\warning\fi
\fi
\let\addcode\gobble % do not write to the closed file again

% The typesetting stage
\typesettingtrue
\def\importlabels#1{} % stop collecting labels and bibliographic items from other documents
\ifscroll \vsize\maxdimen\inewbox\abox\output{\setbox\abox\vbox{\unvbox255\unskip}\shipbox\abox}%
\else % real typesetting
	\ifx\fmtname\plainfmtname
	\totalht\plaintextheight
	\advance\totalht.1in
	\advance\totalht2\vmargin
	\totalwd\plaintextwidth
	\advance\totalwd2\hmargin
	\setpapersize\totalwd\totalht
	\hoffset-1in
	\advance\hoffset\hmargin
	\voffset-1in
	\advance\voffset\vmargin
	\fi
	\output{\plainoutput}
\fi
\let\blah\undefined % disallow placeholder references
\def\bib#1\par{} % do not collect bibliographic items at this stage
\def\refs{\raggedright\rightskip0em plus \maxdimen \advance\bibindent1em \everypar{}\the\bibt\ignorespaces} % typeset bibliography at this stage
\def\addverunused#1#2{} % do not collect verification data at this stage
\def\addcont#1#2#3{} % do not touch the table of contents at this stage
\def\prepass#1{} % do not collect labels and bibliography from other files at this stage
\parn0
\parbreftrue % typeset paragraph back references at this stage
\everypar{\numpar}
\let\chapname\undefined \secn0 \backref0 % reset the numbering
\figno0
\expandafter\gobbleinit\input\jobname\relax
\ifhmode\par\fi\vfill\eject % for LaTeX
\ifhmode\par\fi\vfill\supereject
\end

\expandafter\def\csname\string∞\endcsname{\ifmmode\infty\else{\tensy1}\fi} % ∞ also works outside of formulas
\mathcode`:="603A % colon is a punctuation mark, not a binary operation
\mathchardef\colon="303A % \colon is a binary operation, not a punctuation mark
\ifx\documentclass\undefined
\mathchardef\rightrightarrows="2813
\else
\input amssym
\def\smash{\relax % \relax, in case this comes first in \halign
  \ifmmode\def\next{\mathpalette\mathsmash}\else\let\next\makesmash
  \fi\next}
\def\makesmash#1{\setbox0\hbox{#1}\finsmash}
\def\mathsmash#1#2{\setbox0\hbox{$\mathsurround0pt #1{#2}$}\finsmash}
\def\finsmash{\ht0=0pt \dp0=0pt \box0 }
\fi

\def\marrow#1{\mathbin{\vcenter{\offinterlineskip\kern5pt \smash{$\rightrightarrows$}#1}}}
\def\rarrow{\kern3.338pt \smash{$\rightrightarrows$}}
\def\doublerightarrow{\rightrightarrows}
\def\triplerightarrow{\marrow{\rarrow}}
\def\quadruplerightarrow{\marrow{\rarrow\rarrow}}
\def\quintuplerightarrow{\marrow{\rarrow\rarrow\rarrow}}

 % support
 % colimit
\def\hocolim{\mathop{\rm hocolim}} % colimit
\def\Sing{\mathop{\rm Sing}} % singular simplices
\def\gs{{\bf\Delta}} % geometric simplex
\def\Z{{\bf Z}} % integers
\def\R{{\bf R}} % reals

\def\Y{{\sf Y}} % Yoneda
\def\PSh{{\sf PSh}} % presheaves
\def\op{{\sf op}} % opposite category
\def\CRing{{\sf CRing}} % commutative rings
\def\Top{{\sf Top}} % topological spaces with numerable covers
\def\Loc{{\sf Loc}} % locales with numerable covers
\def\Man{{\sf Man}} % manifolds
\def\MCart{{\sf MCart}} % multicartesian spaces
\def\Cart{{\sf Cart}} % cartesian spaces
\def\Site{{\sf Site}} % numerable site
\def\sSet{{\sf sSet}} % simplicial sets
\def\Set{{\sf Set}} % sets
\def\Open{{\sf Open}} % category of opens

\def\co{{\cal C}} % shape sheaf

\def\Sd{{\sf Sd}} % barycentric subdivision
\def\Ab{{\sf Ab}} % abelian groups
\def\coCh{{\sf coCh}} % cochain complexes

\def\Cc{{\sf C}^*} % singular cochains
\def\Cont{{\sf C}} % continuous functions
\def\HH{{\sf H}} % cohomology
 % local singular cochains
\def\AS{{\sf AS}} % Alexander–Spanier cochains

 % de Rham differential
\def\id{{\rm id}} % identity

\def\sB{{\sf B}} % B for stacks
\def\tB{{\rm B}} % B for spaces
\def\tW{{\rm \bar W}} % delooping of simplicial groups

\def\cV{{\cal V}} % variety of algebras
\def\cE{{\cal E}} % little cubes operad

\def\Map{\mathop{\rm Map}}
\def\RMap{\mathop{{\bf R}{\rm Map}}}

\def\sp{{\sf C}} % an arbitrary Grothendieck topology on topological spaces or locales
\def\Uf{{\eucalbf U}} % numerable site

\yearkeytrue
%\suppressunusedbibtrue

\bib\CU
John~W.~Tukey.
Convergence and Uniformity in Topology.
Annals of Mathematics Studies 2 (\y{1940}).
\doi:10.1515/9781400882199.

\bib\SHT
Samuel Eilenberg.
Singular homology theory.
Annals of Mathematics (Second Series) 45:3 (\y{1944}), 407--447.
\doi:10.2307/1969185.

\bib\GEC
Jean Dieudonné.
Une généralisation des espaces compacts.
Journal de Mathématiques Pures et Appliquées (Neuvième Série) 23 (\y{1944}), 65--76.
\https://gallica.bnf.fr/ark:/12148/bpt6k9704240z/f79.item.

\bib\StonePPS
Arthur~H.~Stone.
Paracompactness and product spaces.
Bulletin of the American Mathematical Society 54:10 (\y{1948}), 977--983.
\doi:10.1090/s0002-9904-1948-09118-2.

\bib\deRham
André Weil.
Sur les théorèmes de de Rham.
Commentarii Mathematici Helvetici 26 (\y{1952}), 119--145.
\doi:10.1007/bf02564296, \eudml:139040.

\bib\CW
Hiroshi Miyazaki.
The paracompactness of CW-complexes.
Tôhoku Mathematical Journal (Second Series) 4:3 (\y{1952}), 309--313.
\doi:10.2748/tmj/1178245380.

\bib\Note
Ernest Michael.
A note on paracompact spaces.
Proceedings of the American Mathematical Society 4:5 (\y{1953}), 831--838.
\doi:10.1090/s0002-9939-1953-0056905-8.

\bib\Fiber
Witold Hurewicz.
On the concept of fiber space.
Proceedings of the National Academy of Sciences 41:11 (\y{1955}), 956--961.
\doi:10.1073/pnas.41.11.956.

\bib\UniBun
John Milnor.
Construction of universal bundles, II.
Annals of Mathematics (Second Series) 63:3 (\y{1956}), 430--436.
\doi:10.2307/1970012.

\bib\MoritaPPS
Kiiti Morita.
Paracompactness and product spaces.
Fundamenta Mathematicae 50:3 (\y{1962}), 223--236.
\doi:10.4064/fm-50-3-223-236.

\bib\PU
Albrecht Dold.
Partitions of unity in the theory of fibrations.
Annals of Mathematics (Second Series) 78:2 (\y{1963}), 223--255.
\doi:10.2307/1970341.

\bib\PNSMS
Kiiti Morita.
Products of normal spaces with metric spaces.
Mathematische Annalen 154:4 (\y{1964}), 365--382.
\doi:10.1007/bf01362570.

\bib\AT
Edwin~H.~Spanier.
Algebraic Topology.
McGraw--Hill, \y{1966}.
\doi:10.1007/978-1-4684-9322-1.

\bib\NumCov
John~E.~Derwent.
A note on numerable covers.
Proceedings of the American Mathematical Society 19:5 (\y{1968}), 1130--1132.
\doi:10.1090/s0002-9939-1968-0234461-x.

\bib\LAT
Albrecht Dold.
Lectures on Algebraic Topology.
Die Grundlehren der mathematischen Wissenschaften 200 (\y{1972}), Springer.
Reprint: Classics in Mathematics (1995), Springer.
\doi:10.1007/978-3-662-00756-3.

\bib\CCCI
Werner Greub, Stephen Halperin, Ray Vanstone.
Connections, Curvature, and Cohomology.  Volume~I.  De Rham Cohomology of Manifolds and Vector Bundles.
Pure and Applied Mathematics 47 (\y{1972}).
\doi:10.1016/S0079-8169(08)62849-4.

\bib\AKTGSC
Kenneth~S.~Brown, Stephen~M.~Gersten.
Algebraic K-theory as generalized sheaf cohomology.
Lecture Notes in Mathematics 341 (\y{1973}), 266--292.
\doi:10.1007/bfb0067062.

\bib\TG
Nicolas Bourbaki.
Topologie Générale.  Chapitres 5 à 10.
Actualités Scientifiques et Industrielles 1235 (\y{1974}), Hermann.
\doi:10.1007/978-3-540-34486-5.

\bib\CSAF
J.~Peter May.
Classifying spaces and fibrations.
Memoirs of the American Mathematical Society 1:155 (\y{1975}).
\doi:10.1090/memo/0155.

\bib\DMAD
Gunter Bengel, Pierre Schapira.
Décomposition microlocale analytique des distributions.
Journées Équations aux dérivées partielles (\y{1978}), 1--2.
\doi:10.5802/jedp.176.

\bib\BottTu
Raoul Bott, Loring~W.~Tu.
Differential Forms in Algebraic Topology.
Graduate Texts in Mathematics 82 (\y{1982}), Springer.
\doi:10.1007/978-1-4757-3951-0.

\bib\Shape
Sibe Mardešić, Jack Segal.
Shape Theory.
North-Holland Mathematical Library 26 (\y{1982}).
\doi:10.1016/B978-0-444-86286-0.50007-6.

\bib\SimpPresh[1985]
J.~Frederick Jardine.
Simplicial presheaves.
Journal of Pure and Applied Algebra 47:1 (1987) 35--87.
\doi:10.1016/0022-4049(87)90100-9.

\bib\EM
Takao Hoshina.
Extensions of mappings II.
(Topics in General Topology.)
North-Holland Mathematical Library 41 (\y{1989}), 41--80.
\doi:10.1016/S0924-6509(08)70148-X.

\bib\Sheaf
Glen~E.~Bredon.
Sheaf Theory.  Second Edition.
Graduate Texts in Mathematics 170 (\y{1997}), Springer.
\doi:10.1007/978-1-4612-0647-7.

\bib\HS[1998]
André Hirschowitz, Carlos Simpson.
Descente pour les $n$-champs.
\arXiv:math/9807049v3.

\bib\AHTS
Fabien Morel, Vladimir Voevodsky.
${\bf A}^1$-homotopy theory of schemes.
Publications mathématiques de l'I.H.É.S.\ 90:1 (\y{1999}), 45--143.
\doi:10.1007/bf02698831.

\bib\Elephant
Peter~T.~Johnstone.
Sketches of an Elephant.  A Topos Theory Compendium.  II.
Oxford Logic Guides 44 (\y{2002}).

\bib\DHI[2002]
Daniel Dugger, Sharon Hollander, Daniel~C.~Isaksen.
Hypercovers and simplicial presheaves.
Mathematical Proceedings of the Cambridge Philosophical Society 136:1 (2004), 9--51.
\arXiv:math/0205027v2, \doi:10.1017/S0305004103007175.

\bib\STopoi[2002]
Bertrand Toën, Gabriele Vezzosi.
Segal topoi and stacks over Segal categories.
\arXiv:math/0212330v2.

\bib\BBT
Dale Husemöller, Michael Joachim, Branislav Jurčo, Martin Schottenloher.
Basic Bundle Theory and $K$-Cohomology Invariants.
Lecture Notes in Physics 726 (\y{2008}), Springer.
\doi:10.1007/978-3-540-74956-1,
{\sevenrm
\https://www.mathematik.uni-muenchen.de/~schotten/Texte/978-3-540-74955-4_Book_LNP726.pdf.
}

\bib\Linf[2008]
Hisham Sati, Urs Schreiber, Jim Stasheff.
$L_∞$-algebra connections and applications to String- and Chern-Simons $n$-transport.
(Quantum Field Theory.  Competitive Models.)
Birkhäuser (\y{2009}), 303--424.
\arXiv:0801.3480v2, \doi:10.1007/978-3-7643-8736-5_17.

\bib\CDT[2008]
Vladimir Voevodsky.
Homotopy theory of simplicial sheaves in completely decomposable topologies.
Journal of Pure and Applied Algebra 214:8 (\y{2010}), 1384--1398.
\arXiv:0805.4578v1, \doi:10.1016/j.jpaa.2009.11.004.

\bib\Cech[2010]
Domenico Fiorenza, Urs Schreiber, Jim Stasheff.
Čech cocycles for differential characteristic classes: an ∞-Lie theoretic construction.
Advances in Theoretical and Mathematical Physics 16:1 (\y{2012}), 149--250.
\arXiv:1011.4735v2, \doi:10.4310/atmp.2012.v16.n1.a5.

\bib\FL
Jorge Picado, Aleš Pultr.
Frames and locales.
Frontiers in Mathematics (\y{2012}), Birkhäuser.
\doi:10.1007/978-3-0348-0154-6.

\bib\HTT
Jacob Lurie.
Higher Topos Theory.
April 9, \y{2017}.
\https://www.math.ias.edu/~lurie/papers/HTT.pdf.

\bib\BG
Dmitri Pavlov.
Question 261195 (revision~1) on MathOverflow.
February 2, \y{2017}.
\https://mathoverflow.net/revisions/261195/1.

\bib\CSp[2019]
Daniel Berwick-Evans, Pedro Boavida de Brito, Dmitri Pavlov.
Classifying spaces of infinity-sheaves.
\arXiv:1912.10544v1.

\bib\DC[2021]
Araminta Amabel, Arun Debray, Peter~J.~Haine.
Differential Cohomology: Categories, Characteristic Classes, and Connections.
\arXiv:2109.12250v1.

\bib\EIB[2021]
Hisham Sati, Urs Schreiber.
Equivariant principal ∞-bundles.
\arXiv:2112.13654v3.

\inewtoks\title
\inewtoks\author
\title{Numerable open covers and representability of topological stacks}
\author{Dmitri Pavlov}
\metadata{\the\title}{\the\author}

%$$
{\tabskip0pt plus 1fil \let\par\cr \obeylines \halign to\hsize{\hfil#\hfil
\articletitle \the\title\vadjust{\bigskip}
\chaptertitle \the\author\vadjust{\medskip}
%$$
Department of Mathematics and Statistics, Texas Tech University
\https://dmitripavlov.org/\vadjust{\medskip}
}}

\hyphenation{Grothen-dieck pre-sheaf pre-sheaves}

\abstract Abstract.
We prove that the class of numerable open covers of topological spaces is the smallest class
that contains covers with pairwise disjoint elements and numerable covers with two elements,
closed under composition and coarsening of covers.
We apply this result to establish an analogue of the Brown--Gersten property for numerable open covers of topological spaces:
a simplicial presheaf on the site of topological spaces
satisfies the homotopy descent property for all numerable open covers
if and only if it satisfies it for numerable covers with two elements and covers with pairwise disjoint elements.
We also prove a strengthening of these results for manifolds, ensuring that covers with two elements can be taken to have a specific simple form.
We apply these results to deduce a representability criterion for stacks on topological spaces similar to \arXiv:1912.10544.
We also use these results to establish new simple criteria for chain complexes of sheaves of abelian groups to satisfy the homotopy descent property.
%This article is available at \arXiv:2203.03120v2.

\tsection Contents

\the\cont

\section Introduction

This paper establishes a generation property (^!{generation of numerable open covers}) for ^{numerable open covers} (^!{numerable open covers})
and then applies this property to establish an analogue of the Brown–Gersten property (^!{descent on manifolds}) for topological spaces and smooth manifolds.
In the introduction, \vgentop\ discusses the general topological results of this paper presented in \vmaintheorems,
whereas \vpresheaves\ discusses their application to simplicial presheaves and the theory of homotopy descent presented in \vappone\ and \vapptwo.

\label\gentop
\subsection Generation of numerable open covers

Recall that an open cover of a topological space is {\it numerable\/} (^!{numerable open covers}) if it admits a ^{partition of unity} (^!{partition of unity}).
The main result of the first part of this paper can be stated as follows.

\proclaim Theorem.
^^={generation of numerable open covers}
The class of ^{numerable open covers} (^!{numerable open covers}) of topological spaces (more generally, locales)
is the smallest class~$C$ of open covers that satisfies the following conditions.
\li The class~$C$ contains all ^{numerable open covers} with two elements.
\li The class~$C$ contains all open covers with pairwise disjoint elements.
\li (Coarsening.)
If $U∈C$ and $U$ refines~$V$, then $V∈C$.
Here $U$ and $V$ are open covers of a space~$X$.
\li (Composition.)
If an open cover $\{U_i\}_{i∈I}$ of a space~$X$ belongs to~$C$ and $ι:X→Y$ is an isomorphism, then the open cover $\{ι(U_i)\}_{i∈I}$ of~$Y$ also belongs to~$C$.
If an open cover $\{U_i\}_{i∈I}$ of a space~$X$ belongs to~$C$, and for every $i∈I$ we are given an open cover $\{V_{i,j}\}_{j∈J_i}$ of~$U_i$ that belongs to~$C$,
then the open cover $\{V_{i,j}\}_{i∈I,j∈J_i}$ of~$X$ belongs to~$C$.

\proof Proof.
See \vmainthm, and also \vmainthmlocale\ for a generalization to the case of locales.

Readers familiar with the notion of a ^{site} (^!{site}) will recognize
the coarsening and composition conditions in ^!{generation of numerable open covers} as saturation conditions for a ^{Grothendieck topology} (^!{Grothendieck topology}).
Thus, ^!{generation of numerable open covers} can be reformulated by saying that the class of ^{numerable open covers}
is the smallest ^{Grothendieck topology} on the category of topological spaces that contains all open covers with pairwise disjoint elements
and all ^{numerable open covers} with two elements.

We also prove the following strengthening of ^!{generation of numerable open covers}
for the case of manifolds (topological or smooth),
which are always assumed to be paracompact and Hausdorff, implying that all open covers are ^{numerable}.
More specifically, we restrict to the case of manifolds isomorphic to disjoint unions of cartesian spaces $\R^n$.

\proclaim Theorem.
^^={generation of open covers of manifolds}
The class of open covers of topological spaces homeomorphic to disjoint unions of $\R^n$
is the smallest class~$C$ that satisfies the following conditions.
\li The class~$C$ contains all open covers $\{U,V\}$ of $\R^n$, where
$$U=⋃_{i∈\Z}(4i,4i+3)⨯\R^{n-1},\qquad V=⋃_{i∈\Z}(4i+2,4i+5)⨯\R^{n-1}.$$
\li The class~$C$ contains all open covers with pairwise disjoint elements.
\li The coarsening and composition conditions from ^!{generation of numerable open covers}.
\endlist This statement continues to hold if we use smooth manifolds and smooth maps instead of topological spaces and continuous maps.

\proof Proof.
See ^!{excellent generation of cartesian spaces}.

^!{generation of open covers of manifolds} relies on the following result about Lebesgue numbers of open covers of noncompact manifolds,
which appears to present independent interest and which we were unable to locate in the literature.

\proclaim Proposition.
For any open cover $\{U_i\}_{i∈I}$ of a cartesian space~$\R^n$ we can find a diffeomorphism $f:\R^n→\R^n$
such that the image of~$U$ under~$f$ has Lebesgue number 1 or higher, meaning the $f$-preimage of any open ball of radius~1 is an subset of some~$U_i$ ($i∈I$).
A similar statement holds for smooth or topological manifolds~$M$, where $f:M\to\R^n$ is now a proper embedding.

\proof Proof.
See \vdistortcover.

\label\presheaves
\subsection The Brown–Gersten property for topological spaces and smooth manifolds

A {\it presheaf\/} of sets on a topological space~$X$ is a functor $F:\Open(X)^\op→\Set$, where $\Open(X)$ denotes the poset of open subsets of~$X$.
{\it Sheaves\/} of sets are singled out from presheaves by the {\it descent condition},
which states that the following diagram is a {\it limit diagram\/}:
$$F(V)→∏_i F(U_i)\doublerightarrow ∏_{j,k}F(U_j∩U_k).$$
This simply boils down to saying that the left map is injective and its image coincides with the subset of the middle set where the two right maps are equal.

A {\it simplicial presheaf\/} (Jardine [\SimpPresh]) on a topological space~$X$ is a functor $F:\Open(X)^\op→\sSet$, where $\sSet$ denotes the category of simplicial sets.
Among simplicial presheaves, one can single out {\it ∞-sheaves\/} via the {\it Čech homotopy descent condition\/} (see Dugger–Hollander–Isaksen [\DHI, Appendix~A];
the notion was also explored by Hirschowitz–Simpson [\HS]),
which requires the following diagram to be a {\it homotopy limit diagram\/}:
$$F(V)→∏_i F(U_i)\doublerightarrow ∏_{j,k}F(U_j∩U_k)\triplerightarrow ∏_{l,m,n}F(U_l∩U_m∩U_n)\quadruplerightarrow ∏_{p,q,r,s}F(U_p∩U_q∩U_r∩U_s)\quintuplerightarrow⋯.$$
Here, the symbol $∏$ denotes {\it homotopy products}.
The diagram on the right is a cosimplicial object in simplicial sets, and homotopy limits of such diagrams are also known as {\it homotopy totalizations}.

Although the general condition for being a homotopy limit diagram is difficult to state,
it simplifies a lot in two special cases: when pairwise intersections are empty ($U_i∩U_j=∅$ for $i≠j$)
and when the open cover has only two elements.

In the case when pairwise intersections are empty, the terms that intersect two or more open sets simply disappear
and the homotopy descent condition simply states that the canonical map to the homotopy product
$$F(V)→∏_i F(U_i)$$
is a weak equivalence.
The homotopy product is easy to compute: simply replace every $F(U_i)$ with a weakly equivalent Kan complex,
then take the ordinary product.

In the case when the open cover has two elements, the diagram now reads
$$F(V)→F(U_1)⨯F(U_2)\doublerightarrow F(U_1∩U_2).$$
This diagram is a homotopy limit diagram if and only if the commutative square
$$\sqcd{F(V)&\mapright{}&F(U_1)\cr\mapdown{}\quad&&\quad\mapdown{}\cr F(U_2)&\mapright{}&F(U_1\cap U_2)\cr}$$
is a homotopy cartesian square.
Equivalently, the map to the homotopy fiber product
$$F(V)→F(U_1)⨯^h_{F(U_1∩U_2)}F(U_2)$$
must be a weak equivalence, which is essentially the homotopy coherent analogue of the Mayer–Vietoris property.
After replacing $F(U_1∩U_2)$ with a weakly equivalent Kan complex if necessary,
the homotopy fiber product on the right side can be computed using the path space construction:
$$F(U_1)⨯_{F(U_1∩U_2)}F(U_1∩U_2)^{Δ^1}⨯_{F(U_1∩U_2)}F(U_2).$$

Thus, even though the general homotopy descent condition can be quite difficult to manipulate,
in the special cases considered above it is quite tractable.
Apart from the fact that homotopy pullbacks are easier to manipulate, they also have better formal properties than homotopy totalizations.
For example, filtered homotopy colimits of simplicial presheaves commute with homotopy pullbacks, but not with homotopy totalizations.
Therefore, we are faced with a question: for what sites the above cases suffice to verify the general homotopy descent condition?

In algebraic geometry, this has been known for a long time.
The {\it Brown--Gersten property\/} allows one to verify the homotopy descent property for simplicial presheaves
on Noetherian topological spaces of finite Krull dimension
(the Zariski topology being the prime example)
by verifying it for the empty cover and for covers consisting of two elements,
see Brown–Gersten [\AKTGSC, Theorem~4].
Morel–Voevodsky [\AHTS, Proposition~1.16] extended this result to the Nisnevich topology,
where the role of pairs of open subsets is assumed by Nisnevich squares.
Voevodsky [\CDT] developed an abstract theory of completely decomposable topologies,
which encompasses both cases.

Recently, homotopy descent for the site of smooth manifolds
(and related sites, such as the small site of a fixed smooth manifold or the site of manifolds with open embeddings as morphisms)
has been gaining in importance in such fields as differential geometry and differential topology,
with many examples, such as (higher) bundle gerbes, twists for vector bundles, Lie ∞-groupoids, and Lie ∞-algebroids,
being formulated in this language,
see, for example, the recent work of Sati–Schreiber–Stasheff [\Linf] and Fiorenza–Schreiber–Stasheff [\Cech].

One is naturally led to the question whether the Brown--Gersten property holds for the site of smooth manifolds.
However, the most obvious analogue fails in a very simple situation.
Fix a set~$A$ and
consider the presheaf~$F$ of sets on the site of smooth manifolds
that sends a smooth manifold~$M$
to the set of locally constant functions $M→A$
with finitely many distinct values.
This presheaf satisfies descent for all finite covers.
However, it does not satisfy descent for the cover of a manifold with infinitely many connected components by its connected components.
Indeed, the relevant descent object consists of all locally constant functions $M→A$, which can have infinitely many distinct values.

Our goal is to show that the above situation is the worst thing that can happen:
if we have homotopy descent for covers with two elements
and covers of a disjoint union $\coprod_i U_i$ by its components~$U_i$,
then we have descent for all open covers.
This is the content of the main theorem, which is proved in greater generality as \vmainthm:

\proclaim Theorem.
^^={descent on topological spaces}
A simplicial presheaf on the site of topological spaces
satisfies the homotopy descent property for all ^{numerable open covers} (^!{numerable open covers})
if and only if it satisfies homotopy descent for
covers with pairwise disjoint elements and numerable covers with two elements.

\proof Proof.
By ^!{numerable open covers form a coverage}, the class of covers with pairwise disjoint elements and numerable covers with two elements forms a ^{coverage} (^!{coverage}).
The class of simplicial presheaves that satisfy the homotopy descent condition with respect to a given ^{coverage}
is the same as the class of simplicial presheaves that satisfy the homotopy descent condition with respect to the ^{Grothendieck topology} (^!{Grothendieck topology})
^{generated} (^!{generated}) by this ^{coverage}.
As shown in \vrefinecovers\ and \vcomposecovers, any ^{Grothendieck topology} is closed under the operations of coarsening and composition.
By ^!{generation of numerable open covers}, saturating the class of covers in the statement
under the operations of coarsening and composition produces the class of all ^{numerable open covers},
which completes the proof.

\proclaim Theorem.
^^={descent on manifolds}
A simplicial presheaf on the site of smooth manifolds
satisfies the homotopy descent property for all open covers
if and only if it satisfies descent for
covers with pairwise disjoint elements
and open covers $\{U,V\}$ of~$\R^n$, where
$$U=⋃_{i∈\Z}(4i,4i+3)⨯\R^{n-1},\qquad V=⋃_{i∈\Z}(4i+2,4i+5)⨯\R^{n-1}.$$

\proof Proof.
Same as the proof of ^!{descent on topological spaces}, but using ^!{generation of open covers of manifolds} instead of ^!{generation of numerable open covers}.

An even more general version for {\it ^{numerable sites}\/} (^!{numerable site}) is given in ^!{good generation of numerable sites}.
In particular, it works for the site of all smooth manifolds,
where morphisms can be taken to be either arbitrary smooth maps, submersions, immersions, etale maps, or open embeddings.
Instead of working with the site of all smooth manifolds we can restrict to the small site of some fixed smooth manifold.
In particular, if we take open embeddings as morphisms we recover the usual notion of a sheaf on a manifold.
The large site of a fixed smooth manifold can also be used.
Smooth manifolds can be replaced with piecewise linear manifolds or topological manifolds.

Although the result appears to have a rather classical flavor,
the author was unable to locate any reference in the literature
that contains or readily implies it,
and such a reference appears to be unknown to a nontrivial fraction of the mathematical community [\BG].
Recent books like Amabel–Debray–Haine [\DC, Theorem~3.6.1]
only cite a weaker result (see, for example, Berwick-Evans–Boavida de Brito–Pavlov [\CSp, Theorem~5.1]),
which shows that open covers with two elements and increasing chains of open embeddings generate the Grothendieck topology on smooth manifolds.

As an application, we leverage the existing work of Berwick-Evans–Boavida de Brito–Pavlov [\CSp, Theorem~5.1] to establish the following
representability result for stacks on the site $\Top$ of topological spaces with numerable open covers (\vnumerablespace).
See \vcspacemap\ for the full statement and details.

\proclaim Theorem.
Suppose $$F:\Top^\op→\sSet$$ is a simplicial presheaf on the site~$\Top$.
Define $$\co F:\Top^\op→\sSet,\qquad \co F(X)=\hocolim_{n∈Δ^\op}F(\gs^n⨯X),$$
where $\gs^n∈\Top$ denotes the $n$-simplex as a topological space
and the homotopy colimit is modeled by the diagonal of a bisimplicial set.
If $F$ satisfies the homotopy descent condition, then $\co F$ is representable in the following sense:
the natural map of \vcmap
$$\co F(X)\to \Map(\Sing X,\co F(*))$$
is a weak equivalence
whenever $X$ is homotopy equivalent to a cofibrant topological space,
e.g., $X$ is a CW-complex, a topological manifold, or a polyhedron.
In particular, we have a natural bijection of sets
$$F[X]\to [\Sing X,\co F(*)],$$
where $F[X]$ denotes concordance classes of the sheaf~$F$ over~$X$.

The significance of this result is that it allows to easily establish Brown-type results
about the existence of classifying spaces and representability of stacks.
For example, it can be used to give a very short proof (\vprincipalbundles) that for any topological group~$G$ the space $\tB G$ classifies numerable principal $G$-bundles.
Likewise, for any abelian topological group~$A$, the Eilenberg–MacLane space $\tB^d A$ classifies numerable $A$-banded bundle $(d-1)$-gerbes (\vbundlegerbes).

As another application, we give the following criterion for acyclicity of sheaves of abelian groups, see \vchsheaf\ for the proof.
As shown in \vflabbyfinesoft, the condition in this criterion is strictly more general than the usual notions of flabby, supple, and fine sheaves, yet it is easy to verify in practice.

\proclaim Theorem.
Suppose $\Site$ is a numerable site (^!{numerable site})
and $F:\Site^\op→\coCh$ is a presheaf of cochain complexes
such that in every cochain degree~$n$
the presheaf of abelian groups
$F_n:\Site^\op→\Ab$
is a sheaf
and for any $M∈\Site$ and any covering family $\{U→M,V→M\}$ of~$M$ the map
$$F_n(U)⨯F_n(V)→F_n(U∩V),\qquad (u,v)↦u|_{U∩V}-v|_{U∩V}$$
is surjective.
Then $F$ satisfies the homotopy descent property.

\subsection Acknowledgments

The author originally wrote up a proof of \vmainthm\ in September 2014 for the paper Berwick-Evans–Boavida de Brito–Pavlov [\CSp].
Subsequent polling on MathOverflow [\BG] revealed that the result is not present in the literature.
Eventually, the authors of [\CSp] figured out a way to prove the main theorem of [\CSp] without this intermediate step,
using a weaker lemma [\CSp, Theorem~5.1], which resulted in the original argument being split off as this article.

I thank Daniel Berwick-Evans and Pedro Boavida de Brito for helpful discussions.
I thank Urs Schreiber for a discussion that led to ^!{excellent generation of cartesian spaces} and ^!{deloopings satisfy descent}.
I thank David Michael Roberts for a discussion that led to ^!{sites of Lie groupoids}.
I thank the anonymous referee of {\it Topology and its Applications\/} for helpful suggestions.

\section Recollections on Grothendieck topologies

In this section we very briefly review the necessary definitions and facts about Grothendieck topologies.
By Toën–Vezzosi [\STopoi, Definition~3.3.1] (see also Lurie [\HTT, Remark~6.2.2.3]),
on an ordinary category (as opposed to an ∞-category)
Grothendieck topologies in the sense of ∞-categories coincide with ordinary Grothendieck topologies.
The ∞-categories that we consider are ordinary categories,
so the usual notion of a Grothendieck topology suffices for our purposes.

\proclaim Definition.
(Johnstone [\Elephant, Definition~C.2.1.1].)
A ^={coverage[|s]} on a category~$C$
is an assignment to every object~$X∈C$
of a collection of families of morphisms (known as ^={generating covering famil[ies|y]}) with codomain~$X$
such that for any such ^{generating covering family} $\{f_i: U_i→X\}_{i∈I}$
and any morphism $g:Y→X$
there is a ^{generating covering family} $\{h_j: V_j→Y\}_{j∈J}$
such that for any $j∈J$ the morphism $gh_j$ factors through the morphism $f_i$ for some $i∈I$.
A ^={site[|s]} is a category with a ^{coverage}.

\proclaim Definition.
A ^={sieve[|s]} on an object~$X$ in a category~$C$
is a collection~$S$ of morphisms with codomain~$X$
such that $f∈S$ implies $fg∈S$ whenever $fg$ is defined.

\proclaim Definition.
(Johnstone [\Elephant, Definition~C.2.1.8].)
A ^={Grothendieck topolog[y|ies]} (alias ^={Grothendieck coverage[|s]})
on a category~$C$ is an assignment to every object~$X∈C$
of a collection of ^{sieves} on~$X$ (known as ^={covering sieve[s|]})
that together form a ^{coverage} on~$C$,
and the following two saturation conditions are satisfied:
for any object $X∈C$ the maximal ^{sieve} (comprising all morphisms with codomain~$X$)
is a ^{covering sieve},
and if $R$ is a ^{covering sieve} on~$X$
and $S$ is a ^{sieve} on~$X$
such that $f^*S=\{g\mid fg∈S\}$ is a ^{covering sieve} for any $f∈R$,
then $S$ is also a ^{covering sieve}.

\proclaim Definition.
The Grothendieck topology on~$C$ ^={generated} by a collection of families of morphisms $\{f_{k,i}:U_i→X_k\}_{k∈K,i∈I_k}$
is the smallest Grothendieck topology (given by the intersection of all such topologies)
for which the given families are contained in some ^{covering sieves}.
In particular, we talk about the ^{Grothendieck topology} and ^{covering sieves} of a ^{site}.

\proclaim Definition.
A ^={covering famil[y|ies]} for a ^{site}~$C$ is a family of morphisms $\{f_i:U_i→X\}_{i∈I}$ in~$C$
such that the smallest ^{sieve} that contains this family is a ^{covering sieve} of~$C$.

The advantage of ^{Grothendieck topologies} is that
passing from a ^{Grothendieck topology} to its category of sheaves
defines an injective map from the collection of ^{Grothendieck topologies} on~$C$
to the collection of full subcategories of the category of presheaves on~$C$.
On the other hand, many different ^{coverages} can have the same category of sheaves.
The advantage of ^{coverages} is that
the given ^{generating covering families} can often be easier
to work with than arbitrary ^{covering sieves}.
For example, the ^{coverage} of $\Man$ (smooth manifolds) consisting of good open covers
(covers whose finite intersections are either empty or diffeomorphic to~$\R^n$)
enjoys many special properties not shared by arbitrary ^{covering sieves}.

We spell out explicitly the particular formulation
of the saturation properties for ^{covering families}
that we will be using below.

\label\refinecovers
\proclaim Lemma.
(cf.~Johnstone [\Elephant, Lemma~C.2.1.6(i)].)
Suppose $C$ is a ^{site}.
Given families of morphisms $\{f_i: U_i\to X\}_{i∈I}$ and $\{g_j: V_j\to X\}_{j∈J}$ in~$C$
such that any $f_i$ factors through some~$g_j$,
if the family~$f$ is a ^{covering family}, then so is~$g$.

\label\composecovers
\proclaim Lemma.
(cf.~Johnstone [\Elephant, Lemma~C.2.1.7(i)].)
Suppose $C$ is a ^{site}.
Given a family $\{f_i: U_i\to X\}_{i∈I}$ and an $I$-indexed collection of families $\{g_{i,j}: V_{i,j}\to U_i\}_{j∈J_i}$
of morphisms in~$C$,
if the family~$f$ and all of~$g_i$ are ^{covering families},
then so is the family $\{f_ig_{i,j}: V_{i,j}\to X\}_{i∈I,j∈J_i}$.

\section Numerable covers

In this section we review the definition of the ^{site} $\Top$
of topological spaces, continuous maps, and ^{numerable open covers},
i.e., open covers that admit a ^{subordinate} ^{partition of unity}.
This ^{site} plays a fundamental role in defining bundle-like structures,
such as principal bundles, vector bundles, (higher) bundle gerbes, etc.,
for which we say that a bundle is {\it numerable\/} if it is trivializable over a ^{numerable open cover}.

This is motivated by the fact that we typically expect or require that
such bundle-like structures are classified by a map to some classifying space.
Numerable bundles are stable under base change because ^{numerable covers} are stable under preimages (Tukey [\CU, V.3.7]).
Thus, if the universal bundle is numerable, then any bundle classified by a map to the classifying space (i.e., the base space
of the universal bundle) must necessarily be numerable.

Indeed, the universal bundle over the classifying space of a topological group~$G$
such that $G$ is a countable CW-complex is numerable because
its base space is a CW-complex (Milnor [\UniBun, Theorem~5.2]),
CW-complexes are paracompact Hausdorff spaces (Miyazaki [\CW]),
and any open cover of a paracompact Hausdorff space is ^{numerable} (Michael [\Note, Proposition~2]).
The latter property of paracompact Hausdorff spaces is the reason why paracompact Hausdorff spaces are sometimes used in expositions related to classifying spaces.
However, there is no need for such a restriction if one works with ^{numerable covers}, as we do below.

We start by recalling some properties of ^{partitions of unity} on topological spaces.
It will be convenient to drop the normalization condition $∑_{i∈I}f_i=1$
and instead require that this sum exists and is a strictly positive continuous map.
We refer to such families as ^{positive partitions}.

\proclaim Definition.
A ^={positive partition[|s]} on a topological space~$X$
is a family of continuous real-valued functions $\{f_i:X→[0,∞)\}_{i∈I}$
such that the map $∑_{i∈I}f_i$ is everywhere defined and is a strictly positive continuous map $X→(0,∞)$.

\proclaim Remark.
For any ^{positive partition} $\{f_i:X→[0,∞)\}_{i∈I}$
and for any $I'⊂I$ the sum $∑_{i∈I'}f_i$ is everywhere defined and is a continuous map $X→[0,∞)$.

\proclaim Definition.
The ^={induced open cover} of a ^{positive partition} $\{f_i\}_{i∈I}$ on a topological space~$X$
is defined as the open cover $\{f_i^*(0,∞)\}_{i∈I}$ of~$X$.
We say that a ^{positive partition} $\{f_i:X→[0,∞)\}_{i∈I}$ is ^={compatible} ^^={compatibility}
with an open cover~$\{U_i\}_{i∈I}$ of~$X$ if $f_i^*(0,∞)⊂U_i$ for all $i∈I$
and ^={subordinate} ^^={subordination} to~$U$ if the closure of $f_i^*(0,∞)$ in~$X$ is a subset of~$U_i$ for all $i∈I$.

\proclaim Remark.
The above definition could be extended to the case when the indexing sets of ^{positive partitions} and open covers are different.
For example, we could say that a ^{positive partition} $\{f_i:X→[0,∞)\}_{i∈I}$ is ^{compatible}
with an open cover~$\{U_j\}_{j∈J}$ of~$X$ if there is a map of sets $α:I→J$ such that $f_i^*(0,∞)⊂U_{α(i)}$ for all $i∈I$.
However, in this case we can define a new ^{positive partition} $\{g_j:X→[0,∞)\}_{j∈J}$ by setting $g_j=∑_{i:α(i)=j}f_i$.
The resulting ^{positive partition} has the same indexing set as the open cover $\{U_j\}_{j∈J}$ and is ^{compatible} with this open cover.
This reindexing procedure works for all definitions and constructions below,
e.g., a locally finite open cover could be reindexed by taking unions indexed by $α^*\{j\}$, etc.
Because of this, our constructions use a fixed indexing set for simplicity.

\proclaim Definition.
A ^{positive partition} $\{f_i\}_{i∈I}$ on a topological space~$X$ is ^={local[ly|] finite[|ness]} if its ^{induced open cover} is a locally finite open cover,
i.e., the collection of open subsets $V⊂X$ such that $f_i|_V=0|_V$ for all but finitely many $i∈I$ forms an open cover of~$X$.

\proclaim Remark.
Occasionally, a weaker property is used.
A ^{positive partition} $\{f_i\}_{i∈I}$ on a topological space~$X$
is ^={point-finite} if for any point $x∈X$ there is only finitely many $i∈I$ such that $f_i(x)≠0$,
i.e., any $x∈X$ belongs only to finitely many elements of the ^{induced open cover} of~$f$.
We make no use of this notion below.

\proclaim Definition.
A ^={partition[|s] of unity} is a ^{positive partition} such that $∑_{i∈I}f_i=1$.

Often, the condition of ^{local finiteness} is included in the definition of a ^{partition of unity}.
We will use both versions below, so for us it is important to be able to distinguish between them.

The following proposition is due to Michael~R.~Mather, see Dold [\LAT, Proposition~A.2.8].
It shows that any ^{positive partition}~$f$
can be improved to a ^{locally finite} ^{partition of unity}~$g$
that is ^{subordinate} to the ^{induced open cover} of~$f$.
Dold's statement only claims ^{compatibility}, not ^{subordination} (i.e., does not mention closure), but our stronger version follows from his proof.

\label\positivepartition
\proclaim Proposition.
For any ^{positive partition} $\{f_i:X→[0,∞)\}_{i∈I}$ on a topological space~$X$
we can find a ^{locally finite} ^{partition of unity} $\{g_i:X→[0,∞)\}_{i∈I}$ on~$X$
such that $g$ is ^{subordinate} to the ^{induced open cover} of~$f$.

\proof Proof.
We normalize $f_i$ by dividing it by $∑_{i∈I}f_i$,
which converts $\{f_i\}_{i∈I}$ into a ^{partition of unity}.
Set $μ=\sup f_i:X→(0,1]$ and $g_i=\max(0,2f_i-μ):X→[0,1]$.
By construction, the family $\{g_i\}_{i∈I}$ is a ^{positive partition}.
The collection of open sets $V_J=h_J^*(0,∞)$,
where $h_J=μ/2-1+∑_{i∈J}f_i$ for a finite subset $J⊂I$,
forms an open cover of~$X$ because $∑_{i∈I}f_i=1$.
We claim that the open cover $\{V_J\}_J$ exhibits the ^{local finiteness} of the ^{positive partition} $\{g_i\}_{i∈I}$.
Indeed, for any finite subset $J⊂I$ and for any $i∈I∖J$ we have $$f_i|_{V_J}≤∑_{k∈I∖J}f_k=1-∑_{k∈J}f_k<μ|_{V_J}/2,$$ so $g_i|_{V_J}=0$, as desired.
\endgraf
Since $2f_i-μ<0$ on an open subset of~$X$, we must have $2f_i-μ≥0$ on the closure of $g_i^*(0,∞)$.
Since $μ>0$ everywhere, we have $f_i>0$ on the closure of $g_i^*(0,∞)$,
so $g$ is ^{subordinate} to the ^{induced open cover} of~$f$.
Replacing $g_i$ with $g_i/∑_{i∈I}g_i$,
we obtain the desired ^{locally finite} ^{partition of unity}.

Thus, an open cover of a topological space
admits a ^{subordinate} ^{locally finite} ^{partition of unity}
if and only if it admits a ^{compatible} ^{positive partition}.
Such covers play a fundamental role
in general topology, where they are known as ^{normal covers}
(introduced by John~W.~Tukey [\CU, §V.2]),
and in the theory of classifying spaces, where they are known as ^{numerable covers} (following Dold [\PU, Definitions~2.1]).

\proclaim Definition.
An open cover $\{U_i\}_{i∈I}$ of a topological space~$X$ is
^={numerable} if it admits a ^{compatible} ^{positive partition}.
^^={numerable cover[s|]}
^^={numerable open cover[s|]}
^^={numerability}

Below we will also apply this definition to open covers of a subset~$X$ of some topological space~$Y$,
in which case we pass to the induced topology on~$X$ first.

The following proposition is a combination of results of
Arthur~H.~Stone [\StonePPS],
Ernest Michael [\Note],
and Kiiti Morita [\MoritaPPS, Theorem~1.2].

Recall that a ^={cozero set[|s]} is an (open) subset $U⊂X$ such that $U=f^*(0,∞)$ for some continuous map $f:X→[0,∞)$.
Recall also that an open cover $\{U_i\}_{i∈I}$ is a ^={star refinement[|s]}
of an open cover $\{V_j\}_{j∈J}$ if the open cover $$\Bigl\{⋃\nolimits_{k∈I: U_k∩U_i≠∅}U_k\Bigr\}_{i∈I}$$ refines~$V$.
A~family $\{U_i\}_{i∈I}$ of open sets is ^={discrete}
if there is an open cover $\{V_j\}_{j∈J}$ of~$X$ such that for every $j∈J$ we have $U_i∩V_j≠∅$ for at most one $i∈I$.
(Any discrete family is in particular disjoint.)

\label\normalcover
\proclaim Proposition.
(Ernest Michael, Kiiti Morita, Arthur~H.~Stone, see Morita [\MoritaPPS, Theorem~1.2]; also Derwent [\NumCov].)
For an open cover $\{U_i\}_{i∈I}$ of a topological space~$X$ the following properties are equivalent.
\li $U$ is a ^{numerable cover}, i.e., it admits a ^{compatible} ^{positive partition};
\li (Michael [\Note, Proposition~2]; see \vpositivepartition.) $U$ admits a ^{subordinate} ^{locally finite} ^{partition of unity};
\li (Michael [\Note, Proposition~1], Hoshina [\EM, Theorem 1.1 and~1.2].)
$U$ can be refined by an open cover given by the union of a countable collection of ^{discrete} families of cozero sets.
\li (Michael [\Note, Theorem~1], Morita [\PNSMS, Theorem~1.2].)
$U$ can be refined by an open cover given by the union of a countable collection of locally finite families of cozero sets.
\li (Stone [\StonePPS, Theorems 1 and~2].)
$U$ is a ^={normal cover[|s]},
^^={normal open cover[|s]}
meaning there is a sequence $W_0=U$, $W_1$, $W_2$, … of open covers of~$X$ such that $W_{n+1}$ is a star refinement of~$W_n$ for all $n≥0$;
\li $U$ can be refined by the inverse image of an open cover of~$Y$ under some continuous map $X→Y$, where $Y$ is a metrizable topological space;
\li (Mardešić–Segal [\Shape, Lemma I.6.1].)
$U$ can be refined by the inverse image of an open cover of~$Y$ under some continuous map $X→Y$,
where $Y$ is an absolute neighborhood retract, i.e.,
a metrizable topological space~$Y$
such that any closed embedding $Y→Z$
into a metrizable topological space~$Z$
factors through an open subset $U⊂Z$
such that there is a map $U→Y$
for which the composition $Y→U→Y$ is identity;
\li $U$ can be refined by a locally finite ^{normal open cover};
\li $U$ can be refined by a locally finite open cover consisting of ^{cozero sets};

\proclaim Remark.
We point out some incompatible definitions of normal and numerable covers in the literature.
Hurewicz [\Fiber, §5] defines normal covers as open covers consisting of ^{cozero sets}.
Hurewicz uses this property in conjunction with local finiteness.
Spanier [\AT, the paragraph before Lemma~2.7.10] defines numerable covers
as locally finite normal covers in the sense of Hurewicz.
May [\CSAF, §1] defines numerable covers as locally finite open covers consisting of ^{cozero sets}.

\proclaim Remark.
Normal spaces are precisely the topological spaces for which every finite open cover is normal (Tukey [\CU, Theorem V.4.1]).
Likewise, paracompact (i.e., every open cover has a locally finite refinement) Hausdorff spaces are precisely
the topological spaces that satisfy the $\rm T_1$ axiom and for which every open cover
is normal (Stone [\StonePPS, Theorem~1]),
or, equivalently, admits a ^{locally finite} ^{compatible} ^{partition of unity} (Michael [\Note, Proposition~2]).

\label\disjointnumerable
\proclaim Remark.
^={Disjoint open cover[s|]} are always ^{numerable}: take $f_i$ to be the characteristic function of~$U_i$.

Most of \vnormalcover\ plays no role in what follows, except for the equivalence
between first three properties, which we prove separately as \vpositivepartition\ and \vnumerablecountable.
Our proof of \vnumerablecountable\ follows Husemöller–Joachim–Jurčo–Schottenloher [\BBT, Proposition 7.1.2].
Although \vnumerablecountable\ holds for ^{discrete} families instead of disjoint families (with the same proof), the latter case is sufficient for our purposes.

\label\numerablecountable
\proclaim Proposition.
(Michael [\Note, Proposition~1]; Hoshina [\EM, Theorem 1.1 and~1.2].)
Any ^{numerable open cover}~$\{U_i\}_{i∈I}$ of a topological space~$X$ can be refined by a ^{numerable open cover}
$\{V_j\}_{j∈J}$ for which there is a map of sets $λ:J→K$,
where $K$ is countable and for any $k∈K$
the family $\{V_j\}_{λ(j)=k}$ is disjoint.

\proof Proof.
Suppose $\{f_i\}_{i∈I}$ is a ^{locally finite} ^{partition of unity} ^{subordinate} to~$\{U_i\}_{i∈I}$.
Denote by $J$ the collection of all finite nonempty subsets of~$I$.
Denote by $K$ the set $\{1,2,3,…\}$ and by $λ:J→K$ the map that computes cardinality.
Define an open cover~$\{V_j\}_{j∈J}$ of~$X$ by
$$V_j=\{x∈X\mid f_k(x)>f_l(x) \hbox{ for all $k∈j$ and $l∈I∖j$}\},$$
where $V_j$ is open because for any $x∈V_j$ we can find an open neighborhood~$W$
such that $f_l|_W≠0$ for only finitely many $l∈I∖j$, so the condition $f_k(x)>f_l(x)$ is an open condition.
The family $\{V_j\}_{j∈J}$ is an open cover of~$X$.
Define a ^{positive partition} $\{g_j\}_{j∈J}$ exhibiting the ^{numerability} of~$V$ by
$$g_j(x)=\max\left(0,\min_{k∈j,l∈I∖j}(f_k(x)-f_l(x))\right).$$
Again, the minimum exists and $g_j$ is continuous because $\{f_i\}_{i∈I}$ is ^{locally finite}.
Also, $V_j=g_j^*(0,∞)$ and $V_j∩V_{j'}=∅$ whenever $j∖j'≠∅$ and $j'∖j≠∅$, by construction.
In particular, for a fixed cardinality of~$j$, all $V_j$ are disjoint.
Since $V_j⊂U_k$ for any $k∈j$, the open cover $V$ refines the open cover~$U$,
which completes the proof.

\proclaim Proposition.
(Tukey [\CU, V.3.7].)
^^={numerable open covers form a coverage}
^{Numerable open covers} form a ^{coverage} on the category of topological spaces and continuous maps.

\proof Proof.
Suppose $X$ is a topological space,
$\{U_i\}_{i∈I}$ is a ^{numerable open cover} of~$X$,
and $g:Y→X$ is a continuous map.
It suffices to show that $\{g^*U_i\}_{i∈I}$ is a ^{numerable open cover} of~$Y$.
Indeed,
if $\{f_i:X→[0,∞)\}_{i∈I}$ is a ^{positive partition} on~$X$ such that $f_i^*(0,∞)⊂U_i$ for all $i∈I$,
then $\{f_ig:Y→[0,∞)\}_{i∈I}$ is a ^{positive partition} on~$Y$ such that $(f_ig)^*(0,∞)=g^*(f_i^*(0,∞))⊂g^*U_i$.

\label\numerablespace
\proclaim Definition.
Denote by~$\Top$ the ^{site} of topological spaces and continuous maps,
equipped with the ^{coverage} of ^{numerable open covers}.

We finish this section by defining another site to which our main theorem is applicable,
the site of locales and their morphisms.
A ^={locale} is a poset that has finite infima and arbitrary suprema
such that the map $b↦\inf(a,b)=a∧b$ preserves suprema for any fixed~$a$,
and a ^={morphism of locales} is an order-preserving map of posets in the opposite direction that preserves finite infima and arbitrary suprema.
These properties axiomatize the properties of the poset of open subsets of a topological space
and the induced inverse image map on open subsets associated to a continuous map.
In particular, we have a functor from topological spaces to locales.
This functor becomes fully faithful on ^={sober spaces}, a large class of topological spaces
that includes all Hausdorff spaces.
The book by Picado–Pultr [\FL] is a comprehensive introduction to locales,
see, in particular, Theorem IX.2.3.4 there, which is a pointfree analogue of (a part of) \vnormalcover.
To keep the paper accessible to a larger audience,
we formulate our lemmas and propositions using topological spaces.
However, nothing in the statements or proofs depends on the availability
of points, since all claims are formulated using open sets only.
In particular, the use of points in the proof of \vnumerablecountable\ is merely a notational convenience.

\label\numerablelocale
\proclaim Definition.
Denote by~$\Loc$ the ^{site} of locales and their morphisms,
equipped with the ^{coverage} of ^{numerable open covers}.

\label\maintheorems
\section Main theorem for topological spaces

In this section, we prove \vmainthm.
The proof repeatedly uses \vrefinecovers\ and \vcomposecovers\ to show that more and more ^{numerable covers} are ^{covering families}
in the ^{Grothendieck topology} ^{generated} by ^{disjoint open covers} (indexed by arbitrary sets) and ^{numerable open covers} with two elements.
First, in \vzigzag\ we show that `zigzag' covers, i.e., covers $P_0,P_1,\ldots$
where the only nonempty intersections are $P_k∩P_{k+1}$, are ^{covering families}.
Next, using a trick with ^{partitions of unity} explained in \vrefinement,
we can refine any countable cover of the form $U_0\subset U_1\subset\cdots$ by a zigzag cover,
which allows us to show that such increasing chains of open subsets are ^{covering families} (\vchain).
\vfinitecover\ shows that finite ^{numerable covers} are ^{covering families}, by induction on the size of the cover.
This is then used in \vccover\ to show that arbitrary countable covers are ^{covering families},
whereas \vacover\ uses \vccover\ and \vnumerablecountable\ to show that covers of arbitrary cardinality are ^{covering families}.
Together with \vdisjointnumerable, \vacover\ implies \vmainthm.

\proclaim Definition.
Denote by $\sp$ the category of topological spaces and continuous maps equipped with some arbitrary ^{Grothendieck topology}
(not necessarily the ^{Grothendieck topology} generated by the ^{coverage} of ^{numerable open covers}).
More generally, $\sp$ can denote the category of locales, as in \vnumerablelocale, but equipped with some arbitrary ^{Grothendieck topology},
not necessarily given by ^{numerable open covers}.

The following lemma and the underlying notion of a zigzag cover
are essentially due to Greub–Halperin–Vanstone (see the proof of Proposition~I.II in~[\CCCI]).

\label\zigzag
\proclaim Lemma.
If countable ^{disjoint open covers} and ^{numerable open covers} with two elements are ^{covering families} in~$\sp$,
then so is any {\it zigzag cover},
i.e., a ^{numerable open cover} $\{P_0,P_1,P_2,\ldots\}$ of~$X$ such that $P_k\cap P_l=\emptyset$ whenever $|k-l|>1$.

\proof Proof.
The only nontrivial intersections of elements of~$P$ are $P_k\cap P_{k+1}$ for $k\ge0$.
Take $$A=P_0\cup P_2\cup P_4\cup\cdots,$$ the (disjoint) union of even elements of~$P$,
and $$B=P_1\cup P_3\cup P_5\cup\cdots,$$ the (disjoint) union of odd elements of~$P$.
We have $A\cup B=X$ and $A\cap B$ equals the (disjoint) union of $P_k\cap P_{k+1}$ for all $k\ge0$.
The open cover $\{A,B\}$ of~$X$ is ^{numerable}:
if $f_i:X→[0,1]$ is a ^{positive partition} ^{compatible} with~$P$,
then $∑_{i≥0}f_{2i}$ and $∑_{i≥0}f_{2i+1}$ form a ^{positive partition} ^{compatible} with~$\{A,B\}$.
The ^{numerable open cover} $\{A,B\}$ of~$X$ is a ^{covering family} by assumption.
The ^{disjoint open covers} $\{P_0,P_2,P_4,\ldots\}$ of~$A$ and $\{P_1,P_3,P_5,\ldots\}$ of~$B$ are ^{covering families} by assumption.
By \vcomposecovers, $\{P_0,P_1,P_2,\ldots\}$ is a ^{covering family} of~$X$.

The following lemma constitutes the technical core of the proof.
It shows that any increasing chain of open subsets can be refined by a zigzag cover.

\label\refinement
\proclaim Lemma.
Suppose $U_0⊂U_1⊂U_2⊂\cdots$ is an increasing sequence of open subsets of a topological space~$X$ that forms a ^{numerable open cover} of~$X$.
Then there is a ^{numerable open cover}~$P=\{P_0,P_1,P_2,\ldots\}$ of~$X$ such that $P_i\subset U_i$
and $P_i\cap P_j=\emptyset$ whenever $|i-j|>1$.
(Thus, $P$ is a zigzag cover of~$X$ in the sense of \vzigzag.)

\proof Proof.
Choose a ^{partition of unity} $\{f_i:X→[0,1]\}_{i≥0}$ ^{compatible} with~$U$.
Define a sequence $\{g_i\}_{i≥-2}$ of functions $$g_i: X\to[0,1],\qquad g_i=\sum_{0\le k\le i}f_k.$$
Fix some sequences $\alpha$, $\beta$, and~$γ$ of real numbers such that $$1>\alpha_i>\beta_i>γ_i>\alpha_{i+1}>\beta_{i+1}>0$$ for all~$i\ge-2$.
Now take $$P_i=g_{i-2}^*[0,\beta_{i-2})\cap g_i^*(\alpha_i,1]$$ for $i\ge0$.
(As one can see from the formula, here the boundary~$-2$ used above becomes important.)
Observe that $$P_i\subset g_i^*(\alpha_i,1]\subset g_i^*(0,1]\subset\bigcup_{0\le k\le i}f_k^*(0,1]\subset\bigcup_{0≤k≤i}U_k=U_i.$$
\endgraf
Next we show that $X=\bigcup_{i≥0}P_i$,
for which it suffices to show by induction on~$k$ that $⋃_{0≤i≤k}P_i=g_k^*(α_k,1]$,
since $⋃_{k≥0}g_k^*(α_k,1]=X$.
The claim holds for $k=-2$ and $k=-1$ because the union is empty and $g_k=0$.
Thus, suppose $⋃_{0≤i<k}P_i=g_{k-1}^*(α_{k-1},1]$ for some $k≥0$.
We have to show that $g_{k-1}^*(α_{k-1},1]∪P_k=g_k^*(α_k,1]$.
One inclusion is trivial and the other boils down to
$$g_k^*(α_k,1]⊂g_{k-1}^*(α_{k-1},1]∪g_{k-2}^*[0,β_{k-2}),$$
which we strengthen to
$$g_k^*(α_k,1]⊂g_{k-1}^*(α_{k-1},1]∪g_{k-1}^*[0,β_{k-2})=g_{k-1}^*((α_{k-1},1]∪[0,β_{k-2}))=g_{k-2}^*[0,1]=X.$$
\endgraf
Now we demonstrate that $P_i\cap P_j=\emptyset$ for all $i\ge 0$ and $j\ge i+2$.
Once we restrict to $P_i\cap P_j$ we have the following inequalities: $g_i>\alpha_i$, $g_{i-2}<\beta_{i-2}$, $g_j>\alpha_j$, $g_{j-2}<\beta_{j-2}$.
(Abusing notation, $α_k$ and $β_l$ denote constant functions on $P_i∩P_j$ with indicated values.)
Thus $$\alpha_i<g_i\le g_{j-2}<\beta_{j-2}<\alpha_i$$ must hold on $P_i∩P_j$, which implies $P_i∩P_j=∅$.
\endgraf
To show that the open cover~$P$ of~$X$ is ^{numerable},
we construct a ^{positive partition} that is ^{compatible} with~$P$.
Take $h_i=2^{-i}⋅\max(0,β_{i-2}-g_{i-2})⋅\max(0,g_i-α_i)$.
By construction, $P_i=h_i^*(0,∞)$.
Thus, $∑_{i≥0}h_i$ exists and is a strictly positive continuous map $X→(0,1]$ because $P$ covers~$X$ and $h_i≤2^{-i}$.

\label\chain
\proclaim Lemma.
If countable ^{disjoint open covers} and ^{numerable covers} with two elements are ^{covering families} in~$\sp$,
then so is any ^{numerable open cover} $$\{U_0,U_1,U_2,\ldots\}$$ of~$X$
such that $U_0\subset U_1\subset U_2\subset\cdots$.

\proof Proof.
Using \vrefinement, we construct a ^{numerable open cover} $\{P_0,P_1,P_2,\ldots\}$ of~$X$
such that $P_k\subset U_k$ and $P_k\cap P_l=\emptyset$ whenever $|k-l|>1$, so the only nontrivial intersections are $P_k\cap P_{k+1}$ for $k\ge0$.
We have $P_k\subset U_k$ and $P$ is a ^{covering family} by \vzigzag.
By \vrefinecovers\ $U$ is also a ^{covering family}.

\label\finitecover
\proclaim Lemma.
If the empty cover and ^{numerable covers} with two elements are ^{covering families} in~$\sp$,
then so is any finite ^{numerable open cover}.

\proof Proof.
We prove by induction on~$n$ that ^{numerable open covers} with fewer than $n$~elements are ^{covering families}.
The empty cover (of the empty set) is a ^{covering family} by assumption, which establishes the base of the induction.
Suppose all ^{numerable open covers} with at most $n>0$~elements are ^{covering families}
and $\{U_0,…,U_n\}$ is a ^{numerable open cover} of~$X$ with a ^{positive partition} $\{f_i:X→[0,∞]\mid 0≤i≤n\}$.
Set $A=(∑_{i<n}f_i)^*(0,∞)∩⋃_{i<n}U_i$.
The family $\{U_0∩A,…,U_{n-1}∩A\}$ forms a ^{numerable cover} of~$A$
with a ^{positive partition} consisting of the restrictions of $f_0$, …, $f_{n-1}$ because $∑_{i<n}f_i>0$ on~$A$.
This family has $n$ elements, so it is a ^{covering family} by induction.
The family $\{A,U_n\}$ is a ^{numerable cover} of~$X$
with a ^{compatible} ^{positive partition} $∑_{i<n}f_i$ and $f_n$,
hence a ^{covering family} by assumption.
(The support of $∑_{i<n}f_i$ need not be a subset of~$A$, so we only require ^{compatibility}, not ^{subordination}.)
By \vcomposecovers, $\{U_0∩A,…,U_{n-1}∩A,U_n\}$ is a ^{covering family} of~$X$.
Hence $\{U_0,…,U_n\}$ is also a ^{covering family} of~$X$ by \vrefinecovers.

\label\ccover
\proclaim Lemma.
If countable ^{disjoint open covers} and ^{numerable covers} with two elements are ^{covering families} in~$\sp$,
then so is any countable ^{numerable open cover}.

\proof Proof.
Suppose $\{V_i\}_{i≥0}$ is a countable ^{numerable open cover} of a topological space~$X$
with a ^{compatible} ^{positive partition} $\{f_i:X→[0,∞)\}_{i≥0}$.
Set $$U_i=V_0∪V_1∪⋯∪V_i$$ and $$g_i=f_0+⋯+f_i:X→[0,∞)$$ for all $i≥0$.
Set $A_i=g_i^*(0,∞)$.
We have $A_0⊂A_1⊂A_2⊂⋯$ and $⋃_{i≥0}A_i=X$.
Also $\{f_i\}_{i≥0}$ is a ^{positive partition} ^{compatible} with the open cover $\{A_i\}_{i≥0}$,
so $A$ is a ^{numerable open cover} of~$X$.
By \vchain, $A$ is a ^{covering family} of~$X$.
\endgraf
Next, $\{V_0∩A_i,V_1∩A_i,…,V_i∩A_i\}$ is a ^{numerable open cover} of~$A_i$
with a ^{positive partition} given by the restrictions of $f_0$, …, $f_i$.
By \vfinitecover, $\{V_0∩A_i,V_1∩A_i,…,V_i∩A_i\}$ is a ^{covering family} of~$A_i$.
Thus, by \vcomposecovers, $\{V_i∩A_j\}_{0≤i≤j}$ is a ^{covering family} of~$X$.
This ^{covering family} refines the open cover $\{V_i\}_{i≥0}$,
so by \vrefinecovers, $\{V_i\}_{i≥0}$ is a ^{covering family} of~$X$.

\label\acover
\proclaim Lemma.
If ^{disjoint open covers} indexed by arbitrary sets and ^{numerable covers} with two elements are ^{covering families} in~$\sp$,
then so is any ^{numerable open cover}.

\proof Proof.
By \vnumerablecountable,
any ^{numerable open cover}~$\{U_i\}_{i∈I}$ of a topological space~$X$ can be refined by a ^{numerable open cover}
$\{V_j\}_{j∈J}$ for which there is a map of sets $f:J→K$,
where $K$ is countable and for any $k∈K$
the family $\{V_j\}_{f(j)=k}$ is disjoint.
Choose a ^{positive partition} $\{h_j:X→[0,∞)\}_{j∈J}$ ^{compatible} with~$V$.
Set $W_k=⋃_{j:f(j)=k}V_j$ and $h_k=∑_{j:f(j)=k}g_j$.
Then $\{h_k\}_{k∈K}$ is a ^{positive partition} ^{compatible} with~$W$,
so $W$ is a countable ^{numerable open cover}.
By \vccover, $W$ is a ^{covering family} of~$X$.
Furthermore, for every $k∈K$ the family $\{V_j\}_{f(j)=k}$ is a ^{disjoint open cover} of~$W_k$,
so also a ^{covering family} of~$W_k$ by assumption.
Thus, by \vcomposecovers, $\{V_j\}_{j∈J}$ is a ^{covering family} of~$X$ that refines~$U$.
Hence, by \vrefinecovers, $\{U_i\}_{i∈I}$ is also a ^{covering family} of~$X$.

\vacover\ combined with \vdisjointnumerable\ immediately implies our main result.

\label\mainthm
\proclaim Theorem.
The ^{Grothendieck topology} of ^{numerable open covers} on $\Top$ is ^{generated}
by ^{disjoint open covers} (indexed by arbitrary sets) and ^{numerable covers} with two elements.

All proofs in this section were carefully formulated to use only open sets, not individual points of~$X$,
and they work equally well without any modification
in the following generalization of \vmainthm\ to locales.

\label\mainthmlocale
\proclaim Theorem.
The ^{Grothendieck topology} of ^{numerable open covers} on $\Loc$ is ^{generated}
by ^{disjoint open covers} (indexed by arbitrary sets) and ^{numerable covers} with two elements.

The following result is used in the proof of ^!{excellent generation of cartesian spaces} and is of interest on its own.

\label\distortcover
\proclaim Proposition.
Suppose $M$ is a smooth manifold and $\{U_i\}_{i∈I}$ is an open cover of~$M$.
Then $M$ admits a smooth proper embedding~$f$ into a cartesian space~$\R^n$
such that the $f$-preimage of any open ball of radius~1 in~$\R^n$ is an open subset of some~$U_i$.
If $M$ is a cartesian space, the embedding can be chosen to be a diffeomorphism.

\proof Proof.
The case of a general~$M$ reduces to the cartesian case as follows.
Pick an arbitrary smooth proper embedding $f:M→\R^n$.
Consider all open subsets $V⊂\R^n$ whose $f$-preimage is a subset of some~$U_i$.
These subsets form an open cover $\{V_j\}_{j∈J}$ of~$\R^n$.
By the cartesian case established below, we get a diffeomorphism $g:\R^n→\R^n$
such that the $g$-preimage of every open ball of radius~1 in~$\R^n$ is an open subset of some $V_j$.
This implies that $gf$ is the desired smooth proper embedding $M→\R^n$.
\ppar
Now assume $M$ is a cartesian manifold.
We construct the diffeomorphism $f:M→M$ in the form $f(v)=a(v)v$, where $a:M→(0,∞)$ is a smooth function
that is radially symmetric, i.e., invariant under the action of the orthogonal group on~$M$.
Thus, $a(v)=b(‖v‖)$,
where $b:[0,∞)→(0,∞)$ is given by $b(x)=a(x⋅e)$ for some fixed unit vector $e∈M$.
Conversely, the function~$a$ constructed from~$b$ is smooth whenever $b$ is smooth and all derivatives of~$b$ at~0 vanish.
The map~$f$ is a diffeomorphism whenever $xb'(x)+b(x)>0$ for all $x≥0$.
\ppar
Suppose $R>0$ and $x∈M$.
We now identify conditions on~$b$ that guarantee that the $f$-preimage of the open ball of radius~1 centered at $f(x)$
is a subset of the open ball of radius $R>0$ centered at~$x$.
This is true whenever for any $w:[0,R]→M$ parametrizing a smooth curve that starts at some point $x∈M$
and $‖w'(t)‖=1$ for all $t∈[0,R]$
we have $1≤\int_0^R a(w(t)) {\rm d}t$.
The latter condition holds whenever for all $y∈[0,∞)$ such that $|y-‖x‖|≤R$ we have $b(y)≥1/R$.
\ppar
Next, we identify conditions on~$b$ that guarantee that for any $x∈M$ we can find $R>0$
with properties as above such that the open ball of radius~$R$ centered at~$x$ is a subset of some~$U_k$.
Denote by $p:[0,∞)→(0,1]$ the continuous function that sends $d∈[0,∞)$
to the maximum $r∈(0,1]$ such that any open ball of radius~$r$ centered at a point of norm at most~$d$ is a subset of some~$U_k$.
Thus, we would like to find a smooth function $b:[0,∞)→(0,∞)$ with vanishing derivatives at~0
such that for all $x≥0$ we have $xb'(x)+b(x)>0$
and for any $x≥0$ there is $R∈(0,p(x)]$ such that for all $y∈[0,∞)∩[x-R,x+R]$ we have $b(y)≥1/R$.
\ppar
Next, we simplify the above conditions by making certain choices.
Assuming $b'(x)≥0$ for all $x≥0$, the second condition is satisfied and the third condition now reads:
for any $x≥0$ there is $R∈(0,p(x)]$ such that $b(x-R)≥1/R$ if $x≥R$ or $b(0)≥1/R$ if $x<R$.
Since $p$ is a continuous nonincreasing function such that $p(0)>0$, we have $y=p(y)$ for some $y≥0$.
Then for all $x≥y$ we have $x≥p(x)≥R$.
Hence, it suffices to ensure that $b(0)≥1/p(y)$ and $b(x-p(x))≥1/p(x)$ for $x≥y$.
The map $G:[y,∞)→[0,∞)$ ($x↦x-p(x)$) is a homeomorphism.
Denote by $g:[0,∞)→[y,∞)$ the inverse homeomorphism.
The last condition on~$b$ now reads: $b(x)≥1/p(g(x))$ for all $x≥0$.
\ppar
Collecting the above conditions,
we would like to find a smooth function $b:[0,∞)→(0,∞)$ with vanishing derivatives at~0
such that
$b(0)≥1/p(y)$ and for all $x≥0$ we have $b'(x)≥0$ and $b(x)≥1/p(g(x))$.
Such a smooth function~$b$ exists: first, we can easily pick a piecewise linear function~$s:[0,∞)→(0,∞)$
with positive slopes that is constant near~0 and satisfies (say) $s(0)≥1+1/p(y)$ and $s(x)≥1+1/p(g(x))$ for all $x≥0$, since the right side is continuous in~$x$.
Then we can make $s$ smooth by using bump functions to smooth out the points where the left and right derivatives are different.

We now present a variant of \vmainthm\ suitable for the cartesian site,
which is relevant for applications.
See, for example, the recent work of Sati–Schreiber–Stasheff [\Linf] and Fiorenza–Schreiber–Stasheff [\Cech].
We cannot quite deduce the theorem below from \vmainthm\ because the newly constructed open sets in \vmainthm\ need not be isomorphic
to cartesian spaces or disjoint unions thereof.

\proclaim Theorem.
^^={excellent generation of cartesian spaces}
Denote by $\MCart$ the full subcategory of $\Top$ on objects that are given by disjoint unions of topological spaces homeomorphic to some~$\R^n$ ($n≥0$).
The ^{Grothendieck topology} of open covers on $\MCart$ is ^{generated}
by ^{disjoint open covers} (indexed by arbitrary sets) and the two-element open cover of $\R^n$ ($n≥0$) by the subsets
$$U=⋃_{i∈\Z}(4i,4i+3)⨯\R^{n-1},\qquad V=⋃_{i∈\Z}(4i+2,4i+5)⨯\R^{n-1}.$$
This statement continuous to hold if we use smooth manifolds and smooth maps instead of topological spaces and continuous maps.

\proof Proof.
Denote by $G$ the ^{Grothendieck topology} generated by the two classes of covers in the statement.
Suppose $\{U_i\}_{i∈I}$ is an open cover of $M∈\MCart$.
Closing $U$ under passage to open subsets produces an open cover $\{V_j\}_{j∈J}$ of~$M$ that refines~$U$.
By \vrefinecovers, if $V∈G$, then $U∈G$.
Thus, we can assume without loss of generality that the open cover $\{U_i\}_{i∈I}$ is closed under passage to open subsets.
Recall also that a Grothendieck topology is closed under isomorphisms, which in our case are given by homeomorphisms, in fact, diffeomorphisms of cartesian spaces.
After applying \vdistortcover\ and further rescaling~$M$, we can assume without loss of generality that every open cube with sides of length~3 or less belongs to the open cover~$U$.
\ppar
We say that an open subset $W⊂M$ is {\it $U$-compatible\/} if the open cover $\{u∈U\mid u⊂W\}$ of~$W$ belongs to~$G$.
Now we show by induction on~$k≥0$ that for any integer $m_{k+1}$, …, $m_n$,
the open subset $A_{k,m}=\R^k⨯∏_{i>k}(2m_i,2m_i+3)⊂M$ is $U$-compatible.
The case $k=0$ was established when we proved that every open cube with sides of length~3 belongs to the open cover~$U$.
Suppose we already proved the statement for all $k<l$.
Taking $k=l$, consider the open cover of $A_{l,m}=\R^l⨯∏_{i>l}(2m_i,2m_i+3)$
by the open sets
$$P_{l,m}=∐_{j∈\Z}\R^l⨯(4j,4j+3)⨯∏_{i>l+1}(2m_i,2m_i+3),\qquad Q_{l,m}=∐_{j∈\Z}\R^l⨯(4j+2,4j+5)⨯∏_{i>l+1}(2m_i,2m_i+3).$$
The open sets $P_{l,m}$ and $Q_{l,m}$ are objects of $\MCart$ (since $(2m_i,2m_i+3)$ is diffeomorphic to~$\R$)
and so is their intersection
$$P_{l,m}∩Q_{l,m}=∐_{j∈\Z}\R^l⨯(4j+2,4j+3)⨯∏_{i>l+1}(2m_i,2m_i+3).$$
Thus, the open cover $\{P_{l,m},Q_{l,m}\}$ of $A_{l,m}$ is a two-element open cover of the type indicated in the statement.
The individual disjoint summands of $P_{l,m}$ and $Q_{l,m}$ are $U$-compatible by the inductive assumption.
Thus, $P_{l,m}$ and $Q_{l,m}$ are $U$-compatible by \vcomposecovers.
Hence, $A_{l,m}$ is $U$-compatible, again by \vcomposecovers.
This completes the inductive step.
Taking $k=n$, we see that $A_{n,m}=\R^n$ is $U$-compatible, which proves that $U∈G$, completing the proof.

\proclaim Remark.
^^={cartesian and multicartesian}
A version of ^!{excellent generation of cartesian spaces} holds for the site $\Cart$ of cartesian spaces given by objects homeomorphic to some $\R^n$,
once we split $U$, $V$, and $U∩V$ into their individual connected components.
We can formulate this as follows: a presheaf~$F$ on $\Cart$ is a sheaf if and only if for every $n≥0$ the commutative square
$$\sqcd{F(\R^n)&\mapright{}&∏_{i∈\Z}F((4i,4i+3)⨯\R^{n-1})\cr
\mapdown{}\quad&&\qquad\mapdown{}\cr
∏_{i∈\Z}F((4i+2,4i+5)⨯\R^{n-1})&\mapright{}&∏_{i∈\Z}F((4i+2,4i+3)⨯\R^{n-1})\cr}$$
is cartesian.
Likewise for smooth manifolds and smooth maps.
This result continues to hold for simplicial presheaves and the homotopy descent condition,
using homotopy products and homotopy cartesian squares in the above statement.

As an illustration of ^!{excellent generation of cartesian spaces} and ^!{cartesian and multicartesian},
we cite the following result of Sati–Schreiber [\EIB, Lemma~3.3.28 and Proposition~3.3.29].
The specific version given below arose from a discussion between Urs Schreiber and the author.

\proclaim Proposition.
^^={deloopings satisfy descent}
For any sheaf~$G$ of groups on the site $\Cart$,
its objectwise delooping groupoid ($M↦\tB G(M)$)
satisfies the homotopy descent property.
More generally, given a presheaf~$G$ of simplicial groups on~$\Cart$,
its objectwise delooping ($M↦\tW G(M)$)
satisfies the homotopy descent property whenever $G$ does.
In the above, groups can be replaced by any variety of algebras that contains the signature and relations of groups, e.g., rings or modules.
Even more generally, simplicial groups can be replaced by algebras over any simplicial algebraic theory that contains groups, e.g., connective spectra
or group-like $\cE_n$-spaces, for any $n≥1$.

\proof Proof.
For the case of sheaves of groups, see Sati–Schreiber [\EIB, Lemma~3.3.28].
The proof boils down to showing that every bundle specified using transition functions
for some open cover of~$\R^n$ is isomorphic to the trivial bundle.
The specific simple form of open sets in ^!{excellent generation of cartesian spaces} and ^!{cartesian and multicartesian} allows for a very simple construction of isomorphism functions, as explained there.
\ppar
For the case of simplicial groups, see Sati–Schreiber [\EIB, Proposition~3.3.29],
which reduces the problem to the previous case by exploiting the fact
that the Čech nerve of an open cover with two elements is 1-skeletal.
\ppar
Finally, for the case of simplicial objects in varieties~$\cV$ of algebras containing groups
(and its homotopy coherent generalization)
it suffices to observe that the forgetful functor from simplicial objects in~$\cV$ to simplicial groups
preserves homotopy limits and homotopy sifted colimits,
and therefore commutes with the homotopy limit over~$Δ$ used in the descent map
as well as the delooping functor $\tB$, which is defined as a homotopy sifted colimit over~$Δ^\op$ of homotopy finite products.

\section Numerable sites

We would like to extend \vmainthm\ to other sites, such as smooth manifolds, possibly equipped with additional geometric structures.
To this end we axiomatize the needed properties in the notion of a ^{numerable site}.

The most general notion of a ^{numerable site} uses the localic site $\Loc$ (\vnumerablelocale).
Readers not familiar with locales can substitute the site $\Top$ (\vnumerablespace) instead of $\Loc$ in ^!{numerable site} without losing any essential examples.

\proclaim Definition.
A ^={numerable site[|s]} is a ^{site}~$S$ equipped with
a functor $\Uf:S→\Loc$
such that $\Uf$ admits cartesian lifts for all open embeddings
and $\Uf$ creates ^{covering families},
meaning the $\Uf$-image of a family~$f$ is a ^{covering family} (^!{covering family}) in~$\Loc$ if and only if $f$ is a ^{covering family} in~$S$.

\label\cartesianlifts
\proclaim Remark.
Unfolding ^!{numerable site}, to specify a ^{numerable site} it suffices to give a functor $\Uf:S→\Loc$
and verify the existence of cartesian lifts for all open embeddings.
The latter condition means that for any object~$\hat Z∈S$
and any open embedding $h:Y→\Uf(\hat Z)$ in $\Loc$ there is a cartesian arrow
$\hat h:\hat Y→\hat Z$ in~$S$ such that $\Uf(\hat h)=h$.
By definition, the arrow $\hat h$ is cartesian whenever for any morphism $\hat f:\hat X→\hat Z$ in~$S$
and any morphism $g:\Uf(\hat X)→Y$ in $\Loc$ such that $hg=\Uf(\hat f)$,
there is a unique morphism $\hat g:\hat X→\hat Y$ such that $\hat h\hat g=\hat f$
and $\Uf(\hat g)=g$.

\proclaim Definition.
A family $\{f_i:U_i→X\}_{i∈I}$ of morphisms in a numerable site~$S$
is ^={disjoint}
^^={disjoint covering famil[y|ies]}
^^={disjoint cover[|s]}
if the family $\{\Uf(f_i)→\Uf(X)\}_{i∈I}$ factors through disjoint open inclusions in $\Loc$.

\proclaim Theorem.
^^={good generation of numerable sites}
If $\Uf:S→\Loc$ is a ^{numerable site},
then the ^{Grothendieck topology} on~$S$
is ^{generated} by ^{covering families} with two elements
and ^{disjoint covering families} indexed by a set.

\proof Proof.
Using \vcartesianlifts, the proof of \vmainthmlocale\ and its dependencies continues to work
for ^{numerable sites}:
all constructions of ^{covering families}
first take the $\Uf$-image in~$\Loc$
and construct the desired open covers in $\Loc$.
Then, the above property allows us to lift the newly constructed open inclusions to~$S$,
and, furthermore, ensure that morphisms whose image factors through the new open inclusion
themselves factor through the constructed lift.

We now give some examples of ^{numerable sites}.
Using \vcartesianlifts, it suffices to give a functor $\Uf:S→\Loc$;
the existence of cartesian lifts of open embedding is obvious for all examples given below.

\proclaim Example.
Normal Hausdorff topological spaces, continuous maps, and locally finite covers form a ^{numerable site}.
Indeed, any locally finite open cover of a normal Hausdorff topological space
admits a ^{subordinate} ^{locally finite} ^{partition of unity} (Bourbaki [\TG, Proposition IX.4.3.3]), hence is a ^{numerable cover}.

\proclaim Example.
Paracompact Hausdorff topological spaces (Dieudonné [\GEC, §1]), continuous maps, and all open covers form a ^{numerable site}.
Indeed, by Stone [\StonePPS, Theorem~1] the class of fully normal (i.e., every open cover is normal) Hausdorff spaces
coincides with the class of paracompact (i.e., every open cover has a locally finite refinement) Hausdorff spaces.

\proclaim Example.
Paracompact Hausdorff smooth, PL, or topological manifolds
together with smooth (respectively PL or continuous) maps,
and all open covers form a ^{numerable site}.

\proclaim Example.
For a ^{numerable site}~$S$ and some object $B∈S$ the slice category $S/B$
inherits a ^{Grothendieck topology} from~$S$ via the forgetful functor $S/B→S$.
The resulting ^{site} is known as the ^={big site[|s]} (alias ^={gros site[|s]}) of~$B$.
It is ^{numerable} because the forgetful functor $S/B→S$ creates ^{covering families} by construction,
the functor $\Uf:S→\Loc$ creates ^{covering families} by definition,
hence the composition $S/B→S→\Loc$ also creates ^{covering families}.
Cartesian lifts for $S/B→\Loc$ can be constructed by lifting to~$S$, which automatically provides the desired map to~$B$.

\proclaim Example.
For a ^{numerable site}~$S$ and some object $B∈S$
the ^={little site[|s]} (alias ^={petit site[|s]}) of~$B$
is the subcategory of the slice category $S/B$
consisting of objects given by morphisms in~$S$ whose image under $\Uf:S→\Loc$ is an open embedding
and morphisms are commutative triangles in~$S$ whose image under $\Uf$ is a commutative triangle of open embeddings.
This subcategory is equipped with the induced Grothendieck topology from $S/B$.
This is a ^{numerable site} since cartesian lifts of open inclusions by definition are morphisms in the little site.

\proclaim Example.
^^={sites of presheaves}
^!{big site} and ^!{little site} can be substantially generalized by replacing $B$ with any presheaf of sets on~$S$,
with the slice category $S/B$ interpreted as the comma category with respect to the Yoneda embedding $\Y:S→\PSh(S)$.
Even more generally, we can take $B$ to be any presheaf of groupoids on~$S$, with the slice category $S/B$ interpreted as the Grothendieck construction for~$B$.
Morphisms from $p:\Y(P)→B$ to $q:\Y(Q)→B$ are then given by maps $f:P→Q$ in~$S$ together with an isomorphism
$p→B(f)(q)$ in the groupoid $B(P)$.
For the same reasons as before, the resulting ^{big site} $S/B$ and ^{little site} of~$B$ are ^{numerable sites}.
(Given an object $(q:\Y(Q)→B)∈S/B$ and a morphism $ι:P→Q$ in~$S$ whose $\Uf$-image is an open embedding,
we can construct a cartesian lift of~$ι$ as the object $(p:\Y(P)→\Y(Q)→B)$ equipped with the morphism to $(q:\Y(Q)→B)$ given by the map $ι:P→Q$ in $S$
together with the identity morphism on~$p$ in the groupoid $B(ι)(q)$.)

\proclaim Example.
Specializing to the case when $B$ is the sheaf of numerable principal $G$-bundles on the site~$\Loc$ for some topological group~$G$,
^!{sites of presheaves} shows that ^!{good generation of numerable sites} holds for the site of topological spaces equipped with a principal $G$-bundle.

\proclaim Example.
^^={sites of simplicial presheaves}
^!{sites of presheaves} generalizes further to {\it $∞$-sites\/} of {\it simplicial\/} presheaves on~$S$,
using the observation of Toën–Vezzosi [\STopoi, Definition~3.3.1] (see also Lurie [\HTT, Remark~6.2.2.3])
that a Grothendieck topology on an ∞-category
can be defined as a ^{Grothendieck topology} on its homotopy category.
^!{good generation of numerable sites} continues to hold for such $∞$-sites, with the obvious adjustments.

We now formulate an analogue of ^!{good generation of numerable sites} for the sites $\MCart$ and $\Cart$ (^!{excellent generation of cartesian spaces}).

\proclaim Definition.
A ^={Euclidean site[|s]} (respectively ^={multi-Euclidean site[|s]}) is a ^{site}~$S$ equipped with
a functor $\Uf:S→\Cart$ (respectively $\MCart$)
such that $\Uf$ admits cartesian lifts for all open embeddings
and $\Uf$ creates ^{covering families},
meaning the $\Uf$-image of a family~$f$ is a ^{covering family} (^!{covering family}) in~$\Cart$ (respectively $\MCart$) if and only if $f$ is a ^{covering family} in~$S$.

\proclaim Theorem.
^^={excellent generation of Euclidean sites}
If $\Uf:S→\MCart$ is a ^{multi-Euclidean site} (^!{multi-Euclidean site}),
then the ^{Grothendieck topology} on~$S$
is ^{generated} by the cartesian lifts of
^{covering families} with two elements
of the form $\{U,V\}$ described in ^!{excellent generation of cartesian spaces},
as well as cartesian lifts of ^{disjoint covering families} indexed by a set.
This also holds for $\Cart$ instead of $\MCart$, with adjustments given in ^!{cartesian and multicartesian}.

\proof Proof.
Same as in ^!{good generation of numerable sites}, but using ^!{excellent generation of cartesian spaces} instead of \vmainthm.

\proclaim Example.
^^={Euclidean sites of presheaves}
^!{sites of presheaves} continues to hold in the Euclidean setting.
Thus, the ^{big site} $S/B$ and the ^{little site} of a presheaf of groupoids~$B$ on $\Cart$ (respectively $\MCart$)
are ^{Euclidean sites} (respectively ^{multi-Euclidean sites}),
as in ^!{Euclidean site}.
^!{sites of simplicial presheaves} with simplicial presheaves also continues to work.

\proclaim Example.
Specializing to the case when $B$ is the sheaf of principal $G$-bundles with connections on the site~$\Cart$ for some Lie group~$G$,
^!{sites of presheaves} shows that ^!{excellent generation of Euclidean sites} holds for the site of cartesian manifolds equipped with a principal $G$-bundle with connection.
We can also take $B$ to be the sheaf of Riemannian metrics, complex or conformal structures, etc.

\proclaim Example.
^^={sites of Lie groupoids}
As a special case of ^!{Euclidean sites of presheaves},
consider a Lie groupoid $$G=(G_0,G_1,s:G_1→G_0,t:G_1→G_0,⋅,e).$$
The ^{big site} $\Cart/G$ of~$G$ (^!{big site}) and the ^{little site} of~$G$ (^!{little site})
are defined by applying the cited constructions to the presheaf of groupoids $\Y(G)$ on~$\Cart$ that is given by applying the Yoneda embedding~$\Y$ to the components of~$G$.
They are ^{Euclidean sites}, so ^!{excellent generation of Euclidean sites} holds for them.

\label\appone
\section Homotopy descent for chain complexes of sheaves of abelian groups

In this section, we use \vmainthm\ to establish a simple sufficient
condition (\vchpresheaf, \vchsheaf) for a chain complex of sheaves of abelian groups to satisfy the homotopy descent property.
The surjectivity condition in \vchsheaf\ below implies
that all of the sheaves $F_n$ are acyclic (\vacyclicsheaf),
so the theorem can be seen as an easy way to show that a sheaf of abelian groups is acyclic.
On the other hand, the required property is strictly weaker than
being flabby or fine (see \vflabbyfinesoft),
so our result does not follow from the classical theory.

\label\chpresheaf
\proclaim Theorem.
Suppose $\Uf:\Site→\Loc$ is a ^{numerable site} (^!{numerable site})
and $F:\Site^\op→\coCh$ is a presheaf of cochain complexes.
Then $F$ satisfies the homotopy descent property whenever
the following conditions are satisfied.
\li In every cochain degree~$n$
the presheaf of abelian groups
$F_n:\Site^\op→\Ab$
sends disjoint unions to products.
\li For any $M∈\Site$ and any covering family $\{U→M,V→M\}$ of~$M$ the map
$$F_n(U)⨯F_n(V)→F_n(U∩V),\qquad (u,v)↦u|_{U∩V}-v|_{U∩V}$$
is surjective in every cochain degree~$n$.
(Here the map $U∩V→U$ is constructed as the cartesian lift of $\Uf(U)∩\Uf(V)→\Uf(U)$ and likewise for~$U∩V→V$, see \vcartesianlifts.)
\li For any $M∈\Site$ and any covering family $\{U→M,V→M\}$ of~$M$ the map
$$F(M)→\ker(F(U)⨯F(V)→F(U∩V))$$ is a quasi-isomorphism.

\proof Proof.
According to \vmainthm, the homotopy descent property for~$F$ is equivalent to the homotopy descent property for disjoint covers and covers with two elements.
The homotopy descent property for disjoint covers holds because
homotopy products of cochain complexes can be computed as ordinary products, and the latter are computed degreewise,
so we conclude by the first property in the statement.
\endgraf
The case of covers with two elements means that for any $M∈\Site$
and a cover $\{U,V\}$ of~$M$
the restriction map
$$F(M)→F(U)⨯^h_{F(U∩V)}F(V)$$ is a quasi-isomorphism.
Computing the homotopy pullback as a homotopy equalizer,
we reduce the problem to showing that in the sequence
$$F(M)→F(U)⨯F(V)→F(U∩V)$$
the left map (given by restrictions to $U$ and~$V$) is the homotopy fiber of the right map
(given by $(u,v)↦(u|_{U∩V}-v|_{U∩V})$).
The right map is surjective by assumption.
Degreewise surjections are fibrations in the projective model structure on cochain complexes,
so the homotopy fiber can be computed as the kernel of the right map.
Thus, it remains to show that the map from $F(M)$ to the kernel of the right map is a quasi-isomorphism,
which holds by assumption.

\label\chsheaf
\proclaim Corollary.
Suppose $\Site$ is a numerable site
and $F:\Site^\op→\coCh$ is a presheaf of cochain complexes
such that in every cochain degree~$n$
the presheaf of abelian groups
$F_n:\Site^\op→\Ab$
is a sheaf
and for any $M∈\Site$ and any covering family $\{U→M,V→M\}$ of~$M$ the map
$$F_n(U)⨯F_n(V)→F_n(U∩V),\qquad (u,v)↦u|_{U∩V}-v|_{U∩V}$$
is surjective.
Then $F$ satisfies the homotopy descent property.

\proof Proof.
We apply \vchpresheaf.
The first property holds because $F_n$ is a sheaf for any~$n$.
The second property holds by assumption.
The third property holds because $F_n$ is a sheaf for any~$n$,
so the map under consideration is an isomorphism.

\label\acyclicsheaf
\proclaim Remark.
In the language of classical abelian sheaf cohomology, \vchsheaf\ says that
a sheaf of abelian groups $F:\Site^\op→\Ab$ is acyclic whenever
for any $M∈\Site$ and any covering family $\{U→M,V→M\}$ of~$M$ the map
$$F_n(U)⨯F_n(V)→F_n(U∩V),\qquad (u,v)↦u|_{U∩V}-v|_{U∩V}$$
is surjective for every~$n$.

The following definition generalizes the notion of a fine sheaf to ^{numerable sites}.
Instead of using the original definition as a sheaf whose endomorphism sheaf is soft
(which would require us to define soft sheaves on numerable sites),
we use an equivalent (for paracompact spaces) characterization in terms of partitions of unity.

\label\finesheaf
\proclaim Definition.
Suppose $\Uf:\Site→\Loc$ is a ^{numerable site} (^!{numerable site}).
A sheaf $$F:\Site^\op→\Ab$$ is {\it fine\/}
^^={fine shea[f|ves]}
if any ^{covering family} $\{U→M,V→M\}$ in~$\Site$
admits a {\it partition of unity}, defined as a pair of endomorphisms of~$F$ (as a sheaf of abelian groups):
$$p_U:F→F,\qquad p_V=\id-p_U:F→F$$ ^{subordinate} to the covering family $\{U→M,V→M\}$.
The subordination condition means that
there is an open inclusion $u'⊂\Uf(M)$
such that $\Uf(U)∪u'=\Uf(M)$
and $p_U$ vanishes when restricted along the cartesian lift $U'→M$ of $u'→\Uf(M)$.
Likewise for $p_V$ and~$V$.

\label\softsheaf
\proclaim Remark.
Suppose $\Uf:\Site→\Loc$ is a ^{numerable site} (^!{numerable site})
and $$F:\Site^\op→\Ab$$ is a sheaf that admits a module structure
over the sheaf of rings $\Cont∘\Uf$, where $$\Cont:\Loc^\op→\CRing$$
sends a locale~$L$ to the commutative ring $\Cont(L,\R)$ of real functions on~$L$.
Then $F$ is a fine sheaf (\vfinesheaf).
Indeed, the endomorphisms $p_U$ and $p_V$ are given by multiplication
by the elements of a partition of unity (^!{partition of unity})
subordinate to the open cover $\{\Uf(U),\Uf(V)\}$ of~$\Uf(M)$.
Analogous observations can be made using smooth functions, or other types of functions,
provided that $\Uf$ factors through smooth manifolds or other appropriate categories.

The following definition is due to Bengel–Schapira [\DMAD].
See there for interesting examples of supple sheaves.

\proclaim Definition.
(Bengel–Schapira [\DMAD].)
Suppose $\Uf:\Site→\Loc$ is a ^{numerable site} (^!{numerable site}).
A sheaf $$F:\Site^\op→\Ab$$ is {\it supple\/}
^^={supple shea[f|ves]}
if for any $U∈\Site$, any $V_1,V_2→U$ whose images are open inclusions,
and any $x∈F(U)$ such that $x|_{V_1∩V_2}=0$,
there are $y_1,y_2∈F(U)$
such that $x=y_1+y_2$,
$y_1|_{V_1}=0$, $y_2|_{V_2}=0$.

We point out an important special case in which the additional condition of \vacyclicsheaf\ is guaranteed to hold.
For topological spaces, this result is well known (see, for example, Bredon [\Sheaf, Theorems II.5.5, II.9.11, II.9.16]),
the novelty here is the generalization to arbitrary ^{numerable sites} and a shorter proof.

% flabby or fine → soft → acyclic
% flabby → supple → soft

\label\flabbyfinesoft
\proclaim Corollary.
Suppose $\Uf:\Site→\Loc$ is a ^{numerable site} (^!{numerable site}).
A sheaf $F:\Site^\op→\Ab$ is acyclic (i.e., $F[0]:\Site^\op→\coCh$ satisfies the homotopy descent property) whenever any of the following conditions hold.
\li $F$ is flabby (French: flasque), i.e., all restriction maps are surjections.
\li $F$ is supple (French: souple) in the sense of ^!{supple sheaf}.
\li $F$ is fine (French: fin) in the sense of \vfinesheaf.
(For example, $F$ admits a module structure like in \vsoftsheaf.)

\proof Proof.
According to \vacyclicsheaf, we have to show that the following property holds:
for any $M∈\Site$ and any covering family $\{U→M,V→M\}$ of~$M$,
the map
$$F(U)⊕F(V)→F(U∩V),\qquad (u,v)↦u|_{U∩V}-v|_{U∩V}$$
is degreewise surjective.
\ppar
For a flabby sheaf this is true because the restriction map $F(U)→F(U∩V)$ is already surjective by definition.
\ppar
For a supple sheaf,
pick an arbitrary section $$s∈F(U∩V).$$
Pick a ^{partition of unity} $f_U,f_V:\Uf(M)→[0,1]$ subordinate to $\{\Uf(U),\Uf(V)\}$.
Set $U'=f_U^*(2/3,1]$ and $V'=f_V^*(2/3,1]$.
Denote by $\hat U'$ and $\hat V'$ the cartesian lifts of $U'$ and $V'$ (\vcartesianlifts).
We have $U'∩V'=∅$ and $\{U',V',\Uf(U)∩\Uf(V)\}$ is an open cover of $\Uf(M)$.
By ^!{supple sheaf}, we can find $y_1,y_2∈F(U∩V)$
such that $s=y_1+y_2$, $y_1|_{\hat U'∩V}=0$, $y_2|_{\hat V'∩U}=0$.
Gluing $y_1$ and the zero section over~$\hat U'$, we get $x_1∈F(U)$.
Likewise, gluing $y_2$ and the zero section over~$\hat V'$, we get $x_2∈F(V)$.
By construction, $x_1|_{U∩V}+x_2|_{U∩V}=y_1+y_2=s$, as desired.
\ppar
For a fine sheaf this is established using a standard argument with a ^{partition of unity}, whose existence is guaranteed by the definition of a fine sheaf.
Throughout this entire proof, we only use those morphisms in $\Site$ whose $\Uf$-images are open inclusions in $\Loc$.
We adopt the convention that capital letters are used for objects of $\Site$,
while small letters are used for their $\Uf$-images in $\Loc$.
Thus, we write $u=\Uf(U)$.
Sometimes, we construct an open inclusion $v→w=\Uf(W)$, in which case its cartesian lift is denoted by $V→W$.
Also, $U∩V→M$ denotes the cartesian lift of $\Uf(U)∩\Uf(V)→\Uf(M)$.
In particular, $\Uf(U∩V)=\Uf(U)∩\Uf(V)$.
\ppar
Using ^!{fine sheaf},
choose a ^{partition of unity} $$p_U:F→F,\qquad p_V=\id-p_U:F→F$$ ^{subordinate} to the covering family $\{U→M,V→M\}$.
The subordination condition means that $p_U$ and $p_V$ are morphisms of sheaves of abelian groups
and there is an open inclusion $u'⊂m$
such that $u∪u'=m$
and $p_U$ vanishes when restricted along the cartesian lift $U'→M$;
a similar condition is imposed on $p_V$ and~$V$.
\ppar
Pick an arbitrary section $$s∈F(U∩V).$$
Denote by $t$ the supremum of all opens~$t→m$ such that the restriction of~$s$ along the cartesian lift of $t∩u∩v→m$ vanishes.
Thus, the cartesian lift $T→M$ can be seen as the complement of the support of~$s$,
and the restriction of $s$ along $T→M$ vanishes.
\ppar
Now construct a section $q_U∈F(U)$ by gluing the zero section over $U∩T$ and the section $p_V(s)∈F(U∩V)$.
These two sections are compatible because their restrictions to $U∩T∩V$ are both~0, by definition of~$T$.
Furthermore, $U∩T→U$ and $U∩V→U$ cover~$U$ because their $\Uf$-images
satisfy $(u∩t)∪(u∩v)=u∩(t∪v)=u∩m=u$,
since $t∪v=m$ by definition of $t$ and $p_V$ (we can take $v'=t$ in the definition of $p_V$).
Symmetrically, construct $q_V∈F(V)$ by gluing $0$ over $V∩T$ and $-p_U(s)$ over $U∩V$.
\ppar
The image of $(q_U,q_V)$ under the map
$$F(U)⊕F(V)→F(U∩V),\qquad (u,v)↦u|_{U∩V}-v|_{U∩V}$$
is $s$ because $q_U$ and $q_V$ restrict to $p_V(s)$ respectively $-p_U(s)$ on $U∩V$ by construction,
and $p_V(s)-(-p_U(s))=p_V(s)+p_U(s)=s$.

\section Example: classical cohomology theories on smooth manifolds

We illustrate the previous section by reproving the classical theorems on the equivalence of de Rham, singular, Alexander–Spanier, and sheaf cohomology on smooth manifolds.
In the formalism of simplicial presheaves,
the other classical cohomology theory, \v Cech cohomology, is equivalent to sheaf cohomology for trivial reasons: the \v Cech nerve of an open cover of a smooth manifold~$M$
is weakly equivalent to~$M$ in the local model structure on simplicial presheaves on the site of smooth manifolds, by definition.

The following theorem is due to Weil [\deRham], see also Bott–Tu [\BottTu, Proposition~8.8].
Recall that $Ω(M)$ denotes the de Rham complex of differential forms on a smooth manifold~$M$.

\proclaim Proposition.
The presheaf of cochain complexes
$$Ω:\Man^\op→\coCh,\qquad M↦Ω(M)$$
satisfies the homotopy descent condition on the site $\Man$.
In particular, the integration map yields a natural weak equivalence
$$Ω(M)→\Cc(M,\R)$$
where $\Cc(M,\R)$ denotes the real smooth singular cochain complex of~$M$.

\proof Proof.
Differential $n$-forms form a sheaf of abelian groups for any~$n≥0$,
which admits a module structure over smooth functions via multiplication,
so the first result follows from \vflabbyfinesoft\ using smooth partitions of unity.
For the second result, observe that both sides satisfy homotopy descent
(see \vsingulardescent\ for the right side)
and are $\R$-invariant, i.e., send projection maps $M\times\R\to M$ to weak equivalences.
Thus, it suffices to show the claim for the case $M=\R^0$, in which case both sides are $\R[0]$ and the map is identity.

Denote by $\Cc(M,A)$ the singular cochain complex with coefficients in an abelian group~$A$ on a smooth manifold~$M$.
Singular $n$-simplices can be taken to be either continuous or smooth maps $\gs^n→M$.

\label\singulardescent
\proclaim Proposition.
For any abelian group~$A$,
the presheaf of cochain complexes
$$\Cc:\Man^\op→\coCh,\qquad M↦\Cc(M,A)$$
satisfies the homotopy descent condition.
In particular, singular cohomology is naturally isomorphic to sheaf cohomology.

\proof Proof.
We invoke \vchpresheaf.
The first property is satisfied because a singular cochain on a disjoint union is a collection of singular cochains on individual terms.
The second property is satisfied because the presheaf $F=\Cc(-,A)$ is flabby:
for any open embedding $W→X$,
the restriction map $\Cc(X,A)→\Cc(W,A)$ is surjective
because a singular cochain on~$W$ can be extended to a singular cochain on~$X$ by zeros.
\ppar
The third property
states that for any $M$ and any cover $\{U,V\}$ of~$M$ the restriction map
$$r:F(M)→\ker(F(U)⨯F(V)→F(U∩V))$$ is a quasi-isomorphism.
In our case, the kernel can be described as the cochain complex $\Cc(\{U,V\},A)$
of singular cochains defined on singular simplices in~$X$ that factor through $U$ or~$V$.
Recall (Eilenberg [\SHT, Chapter~IV]) that the iterated barycentric subdivision construction on singular simplices defines a cochain map $$\Sd:\Cc(\{U,V\},A)→\Cc(M,A)$$
such that $r∘\Sd=\id$ and $\id-\Sd∘r=h∘d+d∘h$ for a cochain homotopy~$h$.
Thus, the restriction map is a quasi-isomorphism, which completes the proof.

Denote by $\AS(M,A)$ the Alexander–Spanier cochain complex with coefficients in an abelian group~$A$ on a topological space~$M$.
In cochain degree~$n$, these are given by germs of $A$-valued functions (not necessarily continuous) on $U^{n+1}$ around the diagonal $U→U^{n+1}$.

\proclaim Proposition.
The presheaf of cochain complexes
$$\AS:\Man^\op→\coCh,\qquad M↦\AS(M,A)$$
satisfies the homotopy descent condition.
In particular, Alexander–Spanier cohomology is naturally isomorphic to sheaf cohomology.

\proof Proof.
The conditions of ^!{fine sheaf} are satisfied:
given an open cover $\{U,V\}$ of~$M$,
pick a ^{partition of unity} $f,g:M→[0,1]$ (so that $f+g=1$)
subordinate to $\{U,V\}$,
set $u'=f^*[0,1/3)$ and $v'=g^*[0,1/3)$.
Finally, define $p_U:\AS_n→\AS_n$ by setting $p_U(s)$ to the germ of an $A$-valued function on
$M^{n+1}$
equal to~$s$ on $(f^*(1/2,1])^{n+1}$ and 0 everywhere else.
Likewise for $p_V$, using $g^*$ instead of~$f^*$.
The pullback of $p_U$ to $u'=f^*[0,1/3)$ vanishes and
$U∪u'=M$.
Likewise for $p_V$.
By \vflabbyfinesoft, we deduce that the sheaf $\AS_n$ is acyclic for all~$n$.
Hence, $\AS$ satisfies the homotopy descent condition.

\label\apptwo
\section Application: Classifying spaces for stacks over topological spaces

In this section we apply \vmainthm\ to prove \vcspacemap,
which provides a general criterion for establishing Brown-type representability results
for simplicial presheaves on topological spaces.

\proclaim Definition.
Suppose $$F:\Top^\op→\sSet$$ is a simplicial presheaf on the site~$\Top$.
Define $$\co F:\Top^\op→\sSet,\qquad \co F(X)=\hocolim_{n∈Δ^\op}F(\gs^n⨯X),$$
where $\gs^n∈\Top$ denotes the $n$-simplex as a topological space
and the homotopy colimit is modeled by the diagonal of a bisimplicial set.

The importance of the presheaf $\co F$ lies in the fact
that $\pi_0 \co F(X)$ is precisely the set of {\it concordance classes\/} of sections of~$F$ over~$X$.
Here two sections $a,b∈F(X)$ are {\it concordant\/} if there is a section $c∈F(\R⨯X)$
such that $$c|_{\{0\}⨯X}≃a \qquad {\rm and}\qquad c|_{\{1\}⨯X}≃b,$$ where $≃$ means that the two given vertices of $F(X)$ are in the same connected component.

The following theorem can be seen as the analogue of Berwick-Evans–Boavida de Brito–Pavlov [\CSp, Theorem~1.2]
for arbitrary topological spaces.

\label\cspace
\proclaim Theorem.
Suppose $$F:\Top^\op→\sSet$$ is a simplicial presheaf on the site~$\Top$.
If $F$ satisfies the homotopy descent condition, then so does $\co F$.

\proof Proof.
By \vmainthm, it suffices to show that $\co F$ satisfies the homotopy descent property
for numerable open covers with two elements and for covers with disjoint elements.
These two cases are proved in Berwick-Evans–Boavida de Brito–Pavlov [\CSp, Theorem~5.3 and Lemma~5.6].
Although the results are formulated for manifolds,
the constructions in their proofs rely exclusively on the existence of partitions of unity,
and continue to work unchanged in the setting of numerable open covers of topological spaces.

\label\cmap
\proclaim Definition.
Suppose $$F:\Top^\op→\sSet$$ is a simplicial presheaf on the site~$\Top$.
The natural map
$$\co F(X)→\Map(\Sing X,\co F(*))$$
is defined as follows.
Using the hom-product adjunction, it suffices to define a map
$$\co F(X)⨯\Sing X→\co F(*).$$
Unfolding the definition of~$\co$ as the diagonal of a bisimplicial set, it suffices to define a map of sets natural in~$n∈Δ$:
$$F(\gs^n⨯X)_n⨯X^{\gs^n}→F(\gs^n)_n.$$
Indeed, given a map $\gs^n→X$, we send the given element of $F(\gs^n⨯X)$ to an element of $F(\gs^n⨯\gs^n)$ and then pull back along the diagonal map $\gs^n→\gs^n⨯\gs^n$
to get an element of $F(\gs^n)$.

\proclaim Remark.
Below, we use the notion of a cofibrant topological space in the Serre–Quillen model structure on topological spaces.
All CW-complexes, topological manifolds, and polyhedra are cofibrant topological spaces.

\label\cspacemap
\proclaim Theorem.
Suppose $$F:\Top^\op→\sSet$$ is a simplicial presheaf on the site~$\Top$.
If $F$ satisfies the homotopy descent condition, then $\co F$ is representable in the following sense:
the natural map of \vcmap
$$\co F(X)\to \Map(\Sing X,\co F(*))$$
is a weak equivalence
whenever $X$ is homotopy equivalent to a cofibrant topological space.
Furthermore, if $F$ preserves weak equivalences, then the natural map
$$\co F(X)\to\RMap(\Sing X,\co F(*))$$
is a weak equivalence for an arbitrary topological space~$X$.

\proof Proof.
The statement about $\RMap$ is an immediate consequence of the statement about $\Map$
and the fact that every topological space~$X$ admits a weak equivalence to a cofibrant topological space.
\ppar
By \vcspace, the presheaf $\co F$ satisfies the homotopy descent condition.
Furthermore, $\co F$ is homotopy invariant, as already observed by Morel–Voevodsky (see, for example, Berwick-Evans–Boavida de Brito–Pavlov [\CSp, Corollary~2.4]).
The class of spaces~$X$ for which the natural map
$$\co F(X)\to \Map(X,\co F(*))$$
is a weak equivalence is closed under retracts and transfinite compositions
because weak equivalences of simplicial sets are closed under retracts and transfinite compositions.
It remains to show that for a map $X_{n-1}→X_n$ that is a cobase change of the inclusion $S^{d-1}→D^d$ for some $d≥0$,
if $X_{n-1}$ belongs to the class, then so does~$X_n$.
This follows from the descent property of $\co F$
with respect to the numerable open cover $\{U,V\}$ of $X_n$,
where $V$ is the interior of~$D^d$ inside~$X_n$
and $U$ is the union of $X_{n-1}$ and a small open band around $S^{d-1}$ inside~$D^d$.
The spaces $U$, $V$, and $U∩V$ are homotopy equivalent to $X_{n-1}$, $D^d$, and $S^{d-1}$, which implies the desired result.

\label\principalbundles
\proclaim Example.
Suppose $G$ is a topological group and consider the simplicial presheaf $\sB G$ on the site~$\Top$
that sends a topological space~$X$ to the nerve of the groupoid of numerable principal $G$-bundles over~$X$.
By \vcspacemap, concordance classes of numerable principal $G$-bundles over a cofibrant topological space~$X$ are classified
by the space
$$\hocolim_{n∈Δ^\op} \sB G(\gs^n) ≃ \tB \hocolim_{n∈Δ^\op} G(\gs^n) ≃ \tB \Sing(G),$$
i.e., the classifying space of the topological group~$G$.
Taking $\pi_0$ on both sides,
we recover the classification of numerable principal $G$-bundles up to concordance for an arbitrary topological group~$G$.
As already observed by Steenrod, the notion of concordance for principal $G$-bundles coincides with the notion of isomorphism.
Even better, \vcspacemap\ provides a statement on the level of spaces, not just~$\pi_0$:
the space of maps $$M→\tB \Sing(G)$$ is weakly equivalent to the space $$(\co \sB G)(M)$$ of numerable principal $G$-bundles on~$M$ up to concordance.

\label\bundlegerbes
\proclaim Example.
Suppose $d≥0$ and $A$ is an abelian topological group.
Consider the simplicial presheaf $\sB^d A$ on the site~$\Top$
that sends a topological space~$X$ to the simplicial set of numerable $A$-banded bundle $(d-1)$-gerbes over~$X$.
The latter can be defined most easily as the fibrant replacement (in the Čech-local projective model structure on $\Top$)
of the simplicial presheaf that sends~$X$ to $\tB^d \Cont(X,A)$,
where $\Cont(X,A)$ denotes the discrete abelian group of continuous maps $X→A$.
By \vcspacemap, concordance classes of such gerbes over a cofibrant topological space~$X$ are classified
by the space
$$\hocolim_{n∈Δ^\op} \sB^d A(\gs^n) ≃ \tB^d \hocolim_{n∈Δ^\op} A(\gs^n) ≃ \tB^d \Sing(A),$$
i.e., the $d$-fold delooping of the topological group~$A$.
Taking $\pi_0$ on both sides,
we recover the classification of numerable $A$-banded bundle $(d-1)$-gerbes up to concordance (in this case, isomorphism), for an arbitrary abelian topological group~$A$:
$$\sB^d A[M] ≅ [M, \tB^d \Sing(A)].$$
Even better, \vcspacemap\ provides a statement on the level of spaces, not just~$\pi_0$:
the space of maps $$M→\tB^d \Sing(A)$$ is weakly equivalent to the space $$(\co \sB^d A)(M)$$ of numerable $A$-banded bundle $(d-1)$-gerbes on~$M$ up to concordance.
If $A$ is a discrete abelian group, then $$[M, \tB^d \Sing(A)]≅\HH^d(M,A)$$
is the $d$th singular cohomology group of~$M$ with coefficients in~$A$.
This provides a geometric interpretation for singular cohomology groups of an arbitrary cofibrant topological space.

\section References

\refs